\documentclass{article}

\usepackage{geometry,amssymb,amsthm,amsmath,
graphics,enumerate}
\usepackage{times}
\usepackage{amsfonts}
\usepackage{amscd}
\usepackage{url}
\usepackage{xcolor}

\newcommand\At{\rm At}

\newcommand{\Fr}{Fra{\"i}ss\'{e} }

\newcommand{\gammaover}{\overline \gamma}

  \def\mcolon{\! : \!}

\newfont{\msbm}{msbm10 scaled\magstephalf}
\def\cf{\mathop{\rm cf}}

\newtheorem{theorem}{Theorem}[subsection]
\newtheorem{corollary}[theorem]{Corollary}

\newtheorem{lemma}[theorem]{Lemma}

\newtheorem{definition}[theorem]{Definition}

\newtheorem{claim}[theorem]{Claim}
\newtheorem{fact}[theorem]{Fact}

\newtheorem{remark}[theorem]{Remark}

\newtheorem{observation}[theorem]{Observation}
\newtheorem{notation}[theorem]{Notation}

\newtheorem{example}[theorem]{Example}

\newtheorem{construction}[theorem]{Construction}

\def\Krich{{\bf R}}

\def\b1K{\mbox{\boldmath $K$}_{-1}}
\def\bK{\mbox{\boldmath $K$}}

\def\bB{\mbox{\boldmath $B$}}

\def\Dscr{{\cal D }}
\def\Pscr{{\cal P     }}

\def\bbB{\mathbb{B}}
\def\b2{\bf 2}

\def\Iscr{\Im}

\def\<{\langle}

\def\>{\rangle}

\def\Mbar{\overline{M}}

\def\meet{\wedge}

\def\acl{{\rm acl}}

\def\cl{{\rm cl}}

\newbox\noforkbox \newdimen\forklinewidth
\forklinewidth=0.3pt \setbox0\hbox{$\textstyle\smile$}
\setbox1\hbox to \wd0{\hfil\vrule width \forklinewidth depth-2pt
 height 10pt \hfil}
\wd1=0 cm \setbox\noforkbox\hbox{\lower 2pt\box1\lower
2pt\box0\relax}
\def\unionstick{\mathop{\copy\noforkbox}\limits}
\def\nonfork_#1{\unionstick_{\textstyle #1}}
\setbox0\hbox{$\textstyle\smile$} \setbox1\hbox to \wd0{\hfil{\sl
/\/}\hfil} \setbox2\hbox to \wd0{\hfil\vrule height 10pt depth
-2pt width
               \forklinewidth\hfil}
\newbox\doesforkbox
\setbox\doesforkbox\hbox{\lower 2pt\box1 \lower
2pt\box2\lower2pt\box0\relax}

\newcommand{\NN}{\mbox{\msbm N}}

\def\subm{\prec_{\bK}}
\def\sub'm{\prec_{\bK'}}
\def\grpf #1 #2{{\rm grp}_{#2}(#1)}
\def\spanf #1 #2{{\rm span}_{#2}(#1)}
\def\fldf #1 #2{{\rm fld}_{#2}(#1)}
\def\dclf #1 #2{{\rm dcl}_{#2}(#1)}
\def\rclf #1 #2{{\rm rcl}_{#2}(#1)}
\def\aclf #1 #2{{\rm acl}_{#2}(#1)}
\def\acff #1 #2{{\rm acf}_{#2}(#1)}
\def\strf #1 #2{{\rm str}_{#2}(#1)}
\def\tclf #1 #2{{\rm acf}_{#2}(#1)}
\def\abar{\mbox{\boldmath $a$}}
\def\bbar{{\bf b}}

\def\hbar{{\bf h}}

\def\vbar{{\bf v}}

\def\xbar{{\bf x}}
\def\ybar{{\bf y}}

\def\tp{{\rm tp }}

\date{\today}

\def\mcolon{\!:\!}

\newcommand{\fr}{\mathrm{fr}}

\begin{document}

\title{Hanf numbers for Extendibility and related phenomena}

\author{John T. Baldwin
\\University of Illinois at Chicago\thanks{Research partially supported by Simons
travel grant G5402, G3535.}
\\
Saharon Shelah\thanks{Item 1092 on Shelah's publication list.
Partially supported by European Research Council grant 338821, and by National Science Foundation grant  136974.
}\\Hebrew University of Jerusalem}

%\sidebar{Submitted Oct. 2016; slight revisions in Feb. 2017; more$F^N_n(c) = F^M_n(c)$
%major revisions Mar 2017 after 2nd report.  }

%diverent versions of the paper.}
%
%\sidebar{June 2 incorporating some responses to souldatos}
%
%\sidebar{souldatos May 31, 2016
%
%Some comments:
%1. On the first reading of the paper, I thought on section 4 you will prove something new. But you are just quoting the result from BKL.
%2. On page 2, you mention our work with Koerwien, but the reference is BKS09 (Baldwin-Kolesnikov-Shelah; initials match Baldwin-Koerwien-Souldatos).
%3. In the definition of the Hanf number, "if P (K, λ) holds for some λ ≥ κ then P (K, λ) holds for arbitrarily large λ". You establish something stronger yourselves: if P(K,lambda) is extendibility, then P holds for measureable cardinals. You do not have an if-then statement.
%
%Some questions:
%If kappa is above a measurable, then do we necessarily have a $\aleph_1$-complete ultrafilter on kappa? I can see that the argument 3.0.4 works for measurable cardinals, but if kappa is above a measurable? Or is there a typo on the statement of 3.0.4?
%Is it possible to have a model of ZFC with only one measurable cardinal? Then non-extendibility fails at the measurable cardinal, but what about the rest?
%
%
%It is not clear what the difference is with ahanfmaxwlf2.tex? nor
%which is the most recent update. On July 1,ahanfmaxwlf2.tex was dated
%june 6 and ahanfmaxwlf3.tex june2
%
%  }
\maketitle

\begin{abstract} This paper contains portions of Baldwin's talk at the
Set Theory and Model Theory Conference (Institute for Research in Fundamental Sciences, Tehran, October 2015) and a detailed proof that in a suitable extension of ZFC, there is a complete sentence of $L_{\omega_1,\omega}$ that has maximal models in cardinals cofinal in the first measurable cardinal and, of course, never again.
 \end{abstract}
In this paper we discuss two theorems whose proofs depend on
extensions of the Fra{\" \i}ss\'{e} method.  We prove here the Hanf
number for the property that every model of a (complete) sentence of
$L_{\omega_1,\omega}$ with cardinality $\kappa$ is
extendible\footnote{We say $\bK$ is {\em universally extendible in
$\kappa$} if $M \in \bK$ with $|M| =\kappa$ has a proper
$\subm$-extension in the class. Here, this means has an
$\infty,\omega$-elementary extension.}  is (modulo some mild set theoretic hypotheses
that we  remove in \cite{BaldwinShelahhanfmaxzfc}) the first measurable
cardinal. And we outline the description of an explicit
$L_{\omega_1,\omega}$-sentence $\phi_n$ characterizing $\aleph_n$ for
each $n$. We provide some context for these developments as outlined
in the lectures at IPM\footnote{Due to our tardiness in preparing this paper it could be included in the special volume dedicated to the 2015 conference.}.

The phrase `Fra{\" \i}ss\'{e} construction' has taken many meanings
in the over 60 years since the notion was born \cite{Fraisse} (and
earlier in an unpublished thesis).  There are two major streams. We
focus here on variants in the original construction, which usually use
the standard notion of substructure.  We don't deal here
directly with `Hrushovski constructions' where a specialized notion of
strong submodel  varying with the case plays a central role. An
annotated bibliography of developments of the Hrushovski variant until  2009 appears at \cite{Baldwinfg}.

The first variant we want to consider is the vocabulary. Fra{\"
\i}ss\'{e} worked with a {\em finite, relational} vocabulary.  While
model theory routinely translates between functions and their graphs
and there is usually little distinction between finite and countable
vocabularies; in the infinite vocabulary case such extensions for the Fra{\"
\i}ss\'{e} construction yield
weaker but still very useful consequences.  The second is a
distinction in goal: the construction of {\em complete} sentences of
$L_{\omega_1,\omega}$ (equivalently studying the {\em atomic} models
of a complete first order theory) rather than constructing
$\aleph_0$-categorical theories. This second shift raises new
questions about the cardinality of the resulting models. The
result in Section~\ref{bkl} pins down more precisely the existence spectra for {\em
complete} sentences of $L_{\omega_1,\omega}$.   Section~\ref{mex} expresses
the role of large cardinal axioms in more algebraic terms. Rephrased,
it says that, consistently with the existence of a measurable
cardinal, there is a nicely defined (by a complete sentence of
$L_{\omega_1,\omega}$) class of models that has non-extendible
(maximal) models cofinally below the first measurable. The previous
upper bound for such behavior was $\aleph_{\omega_1}$.

We acknowledge helpful comments by the referee and by  Joel Berman,
Sakai Fuchino, Menachim Magidor,      Ioannis Souldatos and
especially Will Boney.

\section{Hanf numbers and Spectrum functions in infinitary logic}

Recent years have brought a number of investigations of the spectrum
(cardinals in which a propery occurs) for various phenomena and
various sorts of infinitary definable classes. Some of the relevant
phenomena are existence, amalgamation, joint embedding, maximal
models etc. The class might be defined as an abstract
elementary class, the models of a  (complete) sentence of
$L_{\omega_1,\omega}$, etc.

Hanf observed \cite{Hanf} that for any   property $P(\bK,\lambda)$,
where $\bK$  ranges over a {\em set} of classes of models, there is a
cardinal $\kappa =H(P)$ such that $\kappa$ is the least cardinal
satisfying: if $P(\bK,\lambda)$ holds for some $\lambda\geq \kappa$
then $P(\bK,\lambda)$ holds for arbitrarily large $\lambda$. $H(P)$
is called the Hanf number of $P$. e.g.\ $P(\bK,\lambda)$ might be the
property that $\bK$ has a model of power $\lambda$.

%Similarly the
%L\"owenheim number $\ell(P)$ of a set $P$  of %non-empty
%classes is the least cardinal $\mu$ such that any class $\bK \in P$
%that has a model has one of cardinality $\leq \mu$.

 Morley \cite{Morley65a} showed for an
arbitrary sentence of $L_{\omega_1,\omega}(\tau)$ the Hanf number for existence is $\beth_{\omega_1}$
when $\tau$ is countable (More generally, it is
$\beth_{(2^{|\tau|})^+}$. \cite{Shelahbook}); the situation for {\em complete} sentences
is much more complicated.  Knight \cite{Knightex} found the first
complete sentence characterizing $\omega_1$ (i.e. has a model in $\omega_1$ but no larger)  by building on the
construction of many non-isomorphic $\aleph_1$-like linear orderings.
Hjorth found, by a procedure generalizing the \Fr-construction, for
each $\alpha<\omega_1$, a set $S_\alpha$ (finite for finite $\alpha$)
of complete $L_{\omega_1,\omega}$-sentences\footnote{Inductively,
Hjorth shows at each $\alpha$ and each member $\phi$ of $S_\alpha$
one of two sentences, $\chi_\phi, \chi'_\phi$, works as
$\phi_{\alpha+1}$ for $\aleph_{\alpha +1}$.} such that some
$\phi_\alpha\in S_\alpha$ characterizes $\aleph_\alpha$.
%($\phi_\alpha$ has a model of that cardinality but no larger model).
It is conjectured \cite{Souldatoscharpow} that  it may be impossible
to decide in ZFC which sentence works. Baldwin, Koerwien, and
Laskowski \cite{BKL} show a modification of the Laskowski-Shelah
example (see \cite{LaskowskiShelahatom,
 BFKL}) gives a family of $L_{\omega_1,\omega}$-sentences $\phi_r$,
 which
%homogeneously (see Definition~\ref{newdefhom})
characterize
$\aleph_r$ for $r< \omega$.
%That result established for the first
%time in ZFC, the existence of specific
%$L_{\omega_1,\omega}$-sentences $\phi_r$ characterizing $\aleph_r$.
%These last three results all depended on variants of Fra{\"
%\i}ss\'{e}'s.
In Section~\ref{bkl} we sketch the new notion of
n-disjoint amalgamation that plays a central role in \cite{BKL}.

%The {\em finite amalgamation spectrum} of an abstract elementary
%class $\bK$ is the set $X_{\bK}$ of $n<\omega$ (for\ $\aleph_n \geq
%LS(\bK)$),   and $\bK$ satisfies amalgamation\footnote{We say
%amalgamation holds in $\kappa$ in the trivial special case when all
%models in $\kappa$ are maximal.
% %to distinguish that situation
%%from there simply being no models in $\kappa^+$.
%We say amalgamation fails in $\kappa$ if there are no models to
%amalgamate.} in $\aleph_n$.  There are many examples where the finite
%amalgamation spectrum of a complete sentence  of
%$L_{\omega_1,\omega}$ is either $\emptyset$ or $\omega$.
%
%
%
%As detailed in Theorem~\ref{sumup} for each $1\le r < \omega$, we
%gave the first example of such a sentence with a non-trivial
%spectrum: amalgamation holds up to $\aleph_{r-2}$, but fails in
%$\aleph_{r-1}$.  It holds (trivially) in $\aleph_{r}$ (since all
%models are maximal); there is no model in $\aleph_{r+1}$.

   Further  results by \cite{BKS, KLH, BKSoul}, where the hypothesis are
weakened to allow incomplete sentences of $L_{\omega_1,\omega}$ or
even AEC (Abstract Elementary Classes $(\bK,\leq)$ where the
properties of {\em strong substructure}, $\leq$ are defined
axiomatically) are placed in context in \cite{{BBHanf}}.
 Analogous results were proved earlier for {\em incomplete}
sentences  by \cite{BKSoul} who code certain bipartite graphs in way
that determine specific inequalities between the cardinalities of the
two parts of the graph; in this case all models have cardinality less
than $\beth_{\omega_1}$.

All the exotica mentioned here and described in more detail in
\cite{BBHanf} occurs below $\beth_{\omega_1}$.
% we list a number
%of result providing spectra of amalgamation, joint embedding and
%existence of maximal models. All of this  exotica occurs below
%$\beth_{\omega_1}$.
Baldwin and Boney \cite{BBHanf} have shown that the Hanf number for
amalgamation is no more than the first strongly compact cardinal.
This immense gap motivated the current paper. We show that for the
case of universally extendable (every model has a proper extension), there is a smaller gap. There is a {\em complete} sentence of
$L_{\omega_1,\omega}$ which has a maximal model in  cardinals cofinal
in the first measurable (if such exists), but no larger maximal
model.  Is the same true of amalgamation? That is, can amalgamation
eventually behave very differently than it does in small
cardinalities? At the end of this paper we point to the only known
example where amalgamation (for a complete
$L_{\omega_1,\omega}$-sentence) holds on an initial segment then
fails, then holds again; then there are no larger models.

%The immense gap with the results here show how open the amalgamation
%spectra is.  There are three evident areas: a) try to move the
%techniques here beyond $\aleph_\omega$; b) tighten the bounds in
%\cite{BKS,KLH}; c) going beyond $\beth_{\omega_1}$ in ZFC would
%require totally new ideas.

\medskip

%\medskip
%We noted above that if an AEC has disjoint $(\le
%\aleph_{s},2)$-amalgamation it has a model in $\aleph_{s+2}$. Thus,
%on general grounds we knew $\hat{\bK^r}$ fails disjoint $(\le
%\aleph_{r-1},2)$-amalgamation. But to show ordinary $2$- amalgamation
%failed we had to use our particular combinatorics in
%Lemma~\ref{allmax}.2.  We don't have a `soft' argument that
%`ordinary' amalgamation must fail in $\aleph_{r-1}$.  But there is a
%connection between the amalgamation and existence spectra.
%
%
%
%
%
%Consideration of this conjecture for our examples motivated Part 5 of
%Theorem~\ref{sumup}, which with Lemma~\ref{apfew}  gives a second
%proof of Proposition~\ref{atfail}.  We close with two questions.
%
%
%
%\begin{question}\begin{enumerate}
%\item  Is there a (complete) sentence of $L_{\omega_1,\omega}$ which
%    characterizes $\kappa>\aleph_0$ and has fewer than $2^{\kappa}$
%    models of cardinality $\kappa$?
%\item Is there any AEC, in particular  one defined by a complete
%    sentence in $L_{\omega_1,\omega}$, whose finite non-trivial 2-amalgamation
%    spectrum is not an interval?
%
%    Note that we have found a non-interval for 2-amalgamation but
%    only because the only element is the second interval has {\em
%    trivial} amalgamation.
% \end{enumerate}
%\end{question}

\section{Disjoint Amalgamation}

\subsection{Classes determined by finitely generated
structures}\label{wlf}

\setcounter{theorem}{0} The original Fra{\" \i}ss\'{e} construction
took place in a  {\em finite relational} vocabulary and the resulting
infinite structure was $\aleph_0$-categorical for a first order
theory.  We explore here several ways to construct  a countable
atomic model for a first order theory and thus a complete sentence in
$L_{\omega_1,\omega}$.

Recall (e.g.\ chapter 7 of \cite{Baldwincatmon}) that the models of a
complete sentence of $L_{\omega_1,\omega}(\tau)$ are the reducts to
$\tau$ of the atomic (every finite sequence realizes a principal
type)  models of a complete first order theory in a vocabulary
$\tau'$ extending $\tau$. We discuss classes determined by a
countable set of finitely generated models. In Sections~\ref{mex} and
\ref{bkl}, we describe the examples of such classes used to prove our
main results.

\begin{definition}\label{defclass} Fix a countable vocabulary $\tau$ (possibly with function
symbols). Let $(\bK_0,\subseteq)$  denote a {\em countable}
collection of {\em finite} $\tau$-structures and let
$(\widehat{\bK},\subseteq)$ denote the abstract elementary class
containing all structures $M$ such that every {\em finite generated} substructure of $M$ is in $\bK_0$.
\end{definition}

%\sidebar{No -- direct limits}
These classes have syntactic
characterizations.

\begin{lemma}\label{lemclass}
\begin{enumerate}
\item %If $\bK_0$ is closed under substructure
$\widehat{\bK}$ is
    defined by an $L_{\omega_1,\omega}$-sentence $\phi$. \item If
    $\bK_0$ is closed under substructure then $\phi$ may be taken
    universal \cite{Malitzuniv}. \item  $(\bK_0,\subset)$ satisfies
    the axioms for AEC (except for unions under chains.)
    \end{enumerate}
    \end{lemma}

While traditional Fra{\" \i}ss\'{e} classes are closed under
substructure and produce $\aleph_0$-categorical first order
structures, which are {\em uniformly} locally finite, the search for
atomic models  \cite{Hjorth, BFKL, BKL, BSoul} does not always
require closure under substructure and produces a generic structure
which is locally finite but not uniformly so. In Section~\ref{mex},
we expand the subject further by using
    countable collections of {\em finitely generated} rather than finite structures as
    the `Fra{\" \i}ss\'{e} class'.

%\sidebar{really?? We will see that $\subset$ can be replaced with a
%more useful notion of strong submodel later.}

\begin{definition} \label{defgen} Fix a countable vocabulary $\tau$ (possibly with function
symbols). Let $(\bK_0,\le)$  denote a countable collection of finite
$\tau$-structures with  $(\widehat{\bK},\le)$ as in
Definition~\ref{defclass}.

%denote the abstract elementary class containing all structures $M$
%such that $S_{<\aleph_1}(M) \subseteq \bK_0$.
%We additionally assume that $(\bK_0,\le)$ has the joint embedding property
%(JEP).
\begin{enumerate} \item  A model
 $M\in\widehat{\bK}$ is  {\em rich} or  {\em
${\bK_0}$-homogeneous} if for all %finitely generated
$A$ and $B$
in $\bK_0$ with $A\leq B$, every embedding $f:A\rightarrow M$
extends to an embedding $g:B\rightarrow M$. We denote the class
of rich models in $\widehat \bK$ as $\Krich$.

%\sidebar{does this need mod}
\item The model $M\in\widehat{\bK}$ is {\em generic}  if $M$ is
    rich and $M$ is an increasing union of a countable chain of
    finitely generated
substructures, each of which is in $\bK_0$.

\item We let $\Krich$ denote the subclass of $\widehat{\bK}$
    consisting of rich models.
 \end{enumerate}
\end{definition}

In the examples considered here the generic models will always be
countable.

%\sidebar{ In Hodges, the proof is for $\leq$ taken as substructure
%but the argument extends for any notion of strong substructure
%satisfying the relevant properties of an AEC.}

%\begin{definition} Let $(\bK,\leq)$ be an abstract elementary class.
%We say $\bK$ has $(<\lambda,<\mu)$-amalgamation if for any triple
%$(A,B,C)$ of members of $\bK$ with $A\leq B,C$, there is a $D$ in
%$\bK$ and strong embeddings of $B,C$ into $D$ which agree on $A$. The
%amalgamation is {\em disjoint} if the images of $B$ and $C$ intersect
%in the image of $A$.
%\end{definition}

\begin{definition}\label{dap2} An AEC $(\bK,\leq)$ has {$(<\lambda,2)$-disjoint
amalgamation} if for any $A,B,C \in \bK$ with cardinality $< \lambda$
and $A$ strongly embedded in $B,C$, there is a $D$ and strong
embedding of $B,C$ into $D$ that agree on $A$ and such that the
intersection of their ranges is their image of $A$.

$\bK$ has $2$-amalgamation if the ranges of the embedding are allowed to intersect outside of $f(A)$.

$\bK$ has the {\em joint embedding property} (JEP) if any two models can be embedded in some larger $D$.
\end{definition}

%\begin{lemma}  If $\bK$ has joint embedding and $(<\lambda,<\omega)$-disjoint amalgamation
%then  \begin{enumerate}
%\item $\bK$ has a model in $\lambda^+$
%\item $\bK$ has a rich model in $\lambda$.
%\end{lemma}

%Proof.  Let $M_0 = M$ and define $M_{\alpha +1}$ to be the disjoint
%amalgam of $M_\alpha$ with some finite extension

Fra{\" \i}ss\'{e}'s theorem asserted that if a class of finite models in a finite
relational language is closed under substructure  and satisfies AP
and JEP then there is a generic model whose theory is
$\aleph_0$-categorical and quantifier eliminable. The following
extension of Fra{\" \i}ss\'{e}'s theorem is well-known
\cite{Hodgesbook} and the proof is essentially the same.

\begin{lemma}  \label{getrich}  Suppose $\tau$ is countable and $\bK_0$ is a
 countable class of finite {\em or countable} $\tau$-structures   that
 %each
% $\bK_{0,r}$
 satisfies $2$- amalgamation, in particular $(\le \aleph_0, 2)$-disjoint amalgamation, and JEP,
then

\begin{enumerate}
\item A $\bK_0$-generic (and so rich) $\tau$-structure $M$
    exists. \item if $\bK_0$ is closed under substructure, the
    generic is ultra-homogeneous (every isomorphism between
    arbitrary finitely generated substructures extends to an
    automorphism).
  \end{enumerate}
%\sidebar{Is closure under substructure used or desirable?}

% Moreover, if
%$\bK_0$ is separable,
%%represented in a $\bK_0$-generic $M$.
%$M$  is an atomic model of $Th(M)$.  Further, $\Krich={\bf At}$,
%i.e., every rich model $N$ is an atomic model of $Th(M)$.
\end{lemma}

%Proof. The \Fr argument works. $\qed_{\ref{getrich}}$

 A key distinction from the \Fr situation is that in the first order case $\widehat \bK$ doesn't really
 play a role while in the infinitary case it is an important intermediary between the finitely generated
 structures and $\Krich$. \  \Fr
 passes to the first order theory of the generic since it is $\aleph_0$-categorical in {\em first order logic}.
 In our more general situation the generic may be $\aleph_0$-categorical only in $L_{\omega_1,\omega}$.
The Scott sentence of the rich model gives the $L_{\omega_1, \omega}$
sentence we study. As noted at the beginning of this section we may
regard the models as reducts of atomic models of a first order
theory. Thus $\widehat \bK$ may have arbitrarily large models while
$\Krich$ does not; this holds of some
 examples in \cite{Hjorth, BFKL, BKL}.

%\sidebar{Add connections with atomic models.}

\begin{corollary}  \label{buildrich} Suppose $(\bK_0,\le)$  satisfies
the hypotheses of Lemma~\ref{getrich}.
 Fix $\lambda\ge\aleph_0$.
If $\widehat \bK$ has $(\leq\lambda,2)$-amalgamation and has at most
countably many  isomorphism types of countable structures, then
 every $M\in\widehat \bK$ of power $\lambda$ can be extended to a
    rich model $N\in\widehat \bK$, which is also of power $\lambda$.
%\begin{enumerate}
%\item
%%\item and consequently there is a rich model in $\lambda^+$.
%\end{enumerate}
\end{corollary}

Proof.  Given $M\in\widehat \bK$ of power $\lambda$,  construct a
continuous chain $\<M_i:i<\lambda\>$ of elements of $\widehat \bK$,
each of size $\lambda$. At a given stage $i<\lambda$, focus on a
specific finite substructure $A\subseteq M_i$ and  a particular
finite extension $B\in \widehat \bK$ of $A$. If there is an embedding
of $B$ into $M_i$ over $A$, $M_{i+1} = M_i$.  If not, we may assume
$B\cap M_i=A$. Let $M_{i+1}$ be the disjoint amalgamation of $M_i$
and $B$ over $A$. As there are only $\lambda$-possible extensions, we
can, by iterating, organize this construction so that $N=\bigcup
\{M_i:i<\lambda\}$ is rich. $\qed_{\ref{buildrich}}$

\medskip

 Crucially, in
Section~\ref{compsent} the class $\hat \bK$ under consideration will
not satisfy  two-amalgamation even with finite models; but there will
be amalgamation of free structures with finite.

% For 2), iterating this
%procedure $\lambda^+$ times we get a rich model in $\lambda^+$.

%\sidebar{Check this carefully; much more detail in bkl}

\subsection{Atomic Models of First order theories}

We discuss here classes generated by finite (not finitely generated) structures.
Suppose a generic $\tau$-model $M$ exists. When  is $M$  an atomic model of
its first-order $\tau$-theory? As remarked in  Section 2 of \cite{BKL} this
second condition has nothing to do with the choice of  embeddings on
the class $\bK_0$, but rather with the choice of vocabulary. The
following condition is needed when, for some values of $n$, $\bK_0$
has infinitely many isomorphism types of structures of size   $n$ % (generated by $n$-elements).
%We
%extend \cite{BKL} by allowing weakly locally finite models in a
%properly countable vocabulary.

We denote the class of atomic models of a complete first order theory
by ${\bf At}$.

\begin{definition}  \label{separable}  A class $\bK_0$ of finite structures in a countable
vocabulary is {\em separable} if, for each $A\in\bK_0$ and
enumeration $\abar$ of $A$, there is a  quantifier-free first order
formula $\phi_{\abar}(\xbar)$ such that:
\begin{itemize}
\item  $A\models\phi_{\abar}(\abar)$ and
\item  for all $B\in\bK_0$ and all tuples $\bbar$ from $B$,
    $B\models\phi_A(\bbar)$ if and only if $\bbar$ enumerates a
    substructure $B'$ of $B$ and the map $\abar\mapsto\bbar$ is
    an isomorphism.
\end{itemize}
\end{definition}

In practice, we will apply the observation that if for each
$A\in\bK_0$ and enumeration $\abar$ of $A$, there is a
quantifier-free formula $\phi'_{\abar}(\xbar)$ such that there are
only finitely many $B\in \bK_0$ with cardinality $|A|$ that under
some enumeration $\bbar$ satisfy $\phi'_{\abar}(\bbar)$, then $\bK_0$
is separable.

\begin{lemma} \cite{BKL}  \label{crit}  Suppose $\tau$ is countable and $\bK_0$ is a class of finite $\tau$-structures
that is closed under substructure, satisfies amalgamation, and JEP,
then a  $\bK_0$-generic (and so rich) model $M$ exists. Moreover, if
$\bK_0$ is separable,
%represented in a $\bK_0$-generic $M$.
$M$  is an atomic model of $Th(M)$.  Further, $\Krich={\bf At}$,
i.e., every rich model $N$ is an atomic model of $Th(M)$.
\end{lemma}

Proof:   Since the class $\bK_0$  of finite structures is separable
it has countably many isomorphism types, and thus a $\bK_0$-generic
$M$ exists by the usual Fra{\"i}ss{\'e} construction.
%Lemma~\ref{getgen}.
To show that $M$ is an atomic model of $Th(M)$, it suffices to show
that any finite tuple $\abar$ from $M$ can be extended to a larger
finite tuple $\bbar$ whose type is isolated by a complete formula.
Coupled with the fact that $M$ is $\bK_0$-locally finite, we need
only show that for any finite substructure $A\le M$, any enumeration
$\abar$ of $A$ realizes an isolated type. Since every isomorphism of
finite substructures of $M$ extends to an automorphism of $M$, the
formula $\phi_{\abar}(\xbar)$ isolates $\tp(\abar)$ in $M$.

The final sentence follows since any two rich models are
$L_{\infty,\omega}$-equivalent. $\qed_{\ref{crit}}$

\section{Hanf number for All Models Extendible}\label{mex}

We say an abstract elementary class (the models of a complete
sentence in $L_{\omega_1,\omega}$) is {\em universally extendible in
$\kappa$} if every model of cardinality $\kappa$ has a proper strong
extension ($L_{\infty,\omega}$-elementary extension). In this section
we prove the following theorem.

\begin{theorem}\label{getmax} There is a complete sentence $\phi$ of $L_{\omega_1,\omega}$  that
 has arbitrarily large models. But  under reasonable set theoretic
conditions (specified below), we show that for arbitrarily large
$\lambda < \mu$, where $\mu$ is the first measurable cardinal, and
unboundedly many $\lambda$ if there is no measurable cardinal,  $\phi$ has a
maximal  model(with respect to substructure, which in this case means
$\prec_{\infty,\omega}$) with cardinality between $\lambda$ and
$2^\lambda$.
\end{theorem}

%\sidebar{in last statement: (with respect to substructure????}

We remove in \cite{BaldwinShelahhanfmaxzfc} the set theoretic hypotheses by adapting techniques from \cite{Sh309,Gobelshelah} but at the cost of weakening the freeness of the $P_0$-maximal model; see Remark~\ref{context}.

If $|M|$ is at least the first measurable $\mu$, then for any
$\aleph_1$-complete non-principal ultrafilter $\Dscr$ on $\mu$,
$M^{\mu}/\Dscr$ is a proper extension of $M$. This holds because we
can find an $f \in M^{\mu}$ which hits each element $a\in M$ at most
once.  Thus the equivalence class of $f$ cannot be that of any
constant map on $M$ (since $\Dscr$ is non-principal).  On the other
hand, by the {\L}os theorem for $L_{\omega_1,\omega}$, since $\Dscr$
is $\aleph_1$-complete, the ultrapower is a proper
$L_{\omega_1,\omega}$-elementary extension of $M$.
 Thus, we have shown the Hanf number for extendability is at most $\mu$:

\begin{lemma} If $\mu$ is measurable, for any $\phi \in L_{\mu,\mu}$, in particular in $ L_{\omega_1,\omega}$, no model of cardinality
$\geq \mu$ is maximal.
\end{lemma}

The  proof of the converse (Theorem~\ref{getmax}) fills the remainder
of this section. If we only demand the result for an arbitrary
sentence of $L_{\omega_1,\omega}$ there are easy examples. We learned an example
in terms of $\omega$-models (which is easily reinterpreted into
$L_{\omega_1,\omega}$) from Magidor \cite{Magcrm}.
 %Unfortunately, his lecture notes are no longer on line.
 The following sketch of such an    example will suggest some of the key points of the main argument. Note that we write $\Pscr(X)$ for the powerset of $X$.

%\sidebar{get ref}

%\sidebar{SEpt 14-- definitely the example is screwed up.  fixed I think an ohour later.}

\begin{example}\label{incompex} {\rm  Consider a class $\bK$ of
3-sorted structures where:
$P_0$ is a set, $P_1$ is a boolean algebra of subsets of $P_0$ (given by an
    extensional binary $E$) and  $P_2$ is  just a set;
$\{F_n:n<\omega\}$ is a family of unary functions which assigns to each $c\in P_2$, a sequence $F_n(c) \in P_1$. Demand: $\bigwedge_n F_n(c) = F_n(d)$ implies $d=c$.
%  %  There is one crucial further axiom: for every $c$
%%$\bigcap F_n(c) = \emptyset$.
% if a sequence $F_n(c)$ for a fixed $c \in P_3$  has the
%finite intersection property then the intersection is non-empty.
 Let $\psi \in L_{\omega_1,\omega}$ axiomatize $\bK$.
%
%
%
%
%
%
%
%
%
%
%
%
%
%
%%$P_{-1}$ is a copy of $(\omega,<)$,
%% $P_0$ is a set,
%%$P_1$ is a boolean algebra of subsets of $P_0$ (given by
%%    extensional binary $E$,
%%$P_2$ codes a set of countable sequences of elements of  $P_1$, via a
%%        function $f_n(c) = b$ maps $P_3 \times P_2$ into $P_1$ $c \in P_2, b \in
%%        P_1$.)
%
%
%
%% If a sequence $F_n(c)$ for a fixed $c \in P_3$  has the
%%finite intersection property then the intersection is non-empty.
%
%
We claim $M$ is a maximal model of $\mod(\psi)$  with cardinality
$\lambda$ if $\lambda < $ first measurable, % and $\lambda^\omega = \lambda$;
 $|P^M_0| = \lambda$,
$P^M_1 = \Pscr(P^M_0)$, and
 $P^M_2$ codes  each sequence in $ {}^{\omega}(P^M_1)$ via the $F_n$.

Suppose for contradiction that $N$ with $M\precnsim_{\omega_1,\omega} N$ witnesses non-maximality, then the choice of $M$ and the demand imply that there must be an element $a^* \in P^N_0 -P^M_0$.
Then $D =\{b\in P^M_1: E(a^*,b)\}$ is a non-principal
 ultrafilter  on $\lambda = P^M_0$.  To see that $D$ is non-principal, note that if some $b' \in P_1^M$ generated $D$, then $b' \lneq a^*$, contrary (by elementary extension) to their both being atoms.

Since $D$ is $\aleph_1$-incomplete (as $\lambda$ is not measurable) there exists a sequence $\langle b_n:n<\omega\rangle$ of elements of $P^M_1$ with empty intersection.
Since each  countable sequence of subsets of $P^{M}_0$ is coded as $\langle F^M_n(c): n < \omega \rangle$ for some $c\in P^M_2$, there is a $d\in P_2^M$  with $F_n(d) = b_n$ for each $n$. Thus, $M \models \neg(\exists x) \bigwedge E(x, F_n(d))$, while $N\models \bigwedge E(a^*, F_n(d))$.   This contradicts %$L_{\omega_1,\omega}$-elementarity of
$M \prec_{\omega_1,\omega} N$.
%But the witness to incompleteness must be encoded as $F^M_n(c)$ since there is a unique $c$ encoding each countable sequence from $P^M_1$  and then the final axiom guarantees the intersection is realized in $M$.

 There are
 $2^{\aleph_0}$ $2$-types over the empty set, given, for each
 $X\subset \omega$, via
$(c,d)$ realizes $p_X$ iff $X = \{n:F_n(c) \cap F_n(d) \neq \emptyset\}$.
%$c_X$ realizes $p_X$ if and only if $|a \in P^M_0: a \in f(c,n)| \in
%X$.
This implies no sentence satisfied by $M$ can be complete, since a minor variant of Scott's
characterization of countable models shows that a sentence $\psi$ is
complete if and only if only countably many
$L_{\omega_1,\omega}$-types over $\emptyset$ are realized in models
of $\psi$. In Section~\ref{compsentsec} we modify this example to obtain a complete
sentence.}
\end{example}

%; the predicates $P_0, P_1$ will play similar roles; $P_0$
%will be replaced by function $f_n$ and sequences $f_n(c)$ will be
%coded by a $P_2$.

\subsection{Some preliminaries on Boolean Algebras}

There are a number of slightly different jargons among set theorists,
model theorists, category theorists, and Boolean algebraists. In this
section we will spell some of them out, indicate some translations,  specify our notation, and prove some properties of Boolean algebras that will be used in the proof.

An ultrafilter {\em of} a Boolean algebra $B$ is a maximal filter
(i.e. a subset of $B$ that is closed up, under intersection and
contains either $a$ or $a^-$ -- the complement of $a$).  An
ultrafilter {\em on} a set $X$ is a subset of its power set
and  so is an ultrafilter of the Boolean algebra  $\Pscr(X)$.

We begin with some basic properties of independence in Boolean
algebras. A key fact is an equivalence of two notions of independence
on countably infinite Boolean algebras that disappears in the
uncountable. That is, a countable Boolean algebra is
$\aleph_0$-categorical if and only it is free on countably many
generators in the sense of \ref{bn} if and only if it is generated by an
independent set in the  sense of \ref{indepset}. But this equivalence fails in the
uncountable.

\begin{definition}\label{bn}

\begin{enumerate}
\item For $X \subseteq B$ and $B$ a Boolean algebra, $\overline X
    = X_B =\langle X\rangle_B$ be the subalgebra of $B$ generated
    by $X$.
\item A set $Y$ is {\em independent} (or {\em free}) from $X$
   modulo an ideal $\Iscr$ (with domain $I$) in a Boolean algebra $B$ if and only if
    for any Boolean-polynomial $p(v_0,\ldots ,v_k)$ (that is not
    identically $0$), and any
   % \footnote{Shelah April 13, 2017
%    marks out this clause. But I think this is incorrect, without
%    it the equation $x_1 \wedge x_2 =  x_1 \wedge x_2$ would
%    violate the condition even if $X,Y$ were a partition of a set
%    of free generators. But that note also reverses X and Y
%    (actually takes Y as an ideal I) and applies the $\sigma$ to
%    elements of I. None of this makes sense to me.}
   $a\in \langle X\rangle_B- \Iscr$, and distinct $y_i \in Y$,
 $p(y_0,\ldots ,y_k)\wedge a\not \in \Iscr$.

 \item Such an independent $Y$ is called  a {\em basis} for $\langle X \cup Y \cup I\rangle$ over $\langle X \cup I\rangle$.
\end{enumerate}

%A set $Y$ is independent over $X$ in a Boolean algebra if
%if and
%    only if for any Boolean-polynomial $p(v_0,\ldots ,v_k)$,  $B
%    \models p(y_0,\ldots ,y_k)= a$ for some $a \in X$ for some sequence $y_i$ from $Y$
%    implies that $a =0$ and $p(v_0,\ldots ,v_k)=0$ is a law holding in each
%    Boolean algebra.

\end{definition}

There is no requirement that $\Iscr$ be contained in $X$.  Observe the following:

\begin{observation}\label{obs1}
\begin{enumerate}
\item If $\Iscr$ is the $0$ ideal, (i.e., $Y$  is independent from
    $X$), the condition becomes: for any $a\in \langle
    X\rangle_B- \{0\}$, $B \models p(y_0,\ldots ,y_k)\wedge a>
    0$.  That is, every finite Boolean combination of elements of
    $Y$ meets each non-zero $a \in \langle X\rangle_B$.

\item Let $\pi$ map $B$ to $B/\Iscr$. If `Y is independent from
    $X$ over $\Iscr$' then the image of $Y$ is free from the
    image of $X$ (over $\emptyset$) in $B/\Iscr$.  Conversely, if
    $\pi(Y)$ is independent over $\pi(X)$ in $B/\Iscr$, for any
    $Y'$ mapping by $\pi$ to $\pi(Y)$, $Y'$ is independent from
    $X$ over $\Iscr$.

So, if  $X$ is empty, the condition `Y is independent  over
    $\Iscr$' implies  the image of $Y$ is an independent subset
    of $B/\Iscr$.

    \item If a set $Y$ is {\em independent} (or {\em free}) from
        $X$ over an ideal $\Iscr$ in a Boolean algebra $B$ and
        $Y_0$ is a subset of $Y$, then $Y-Y_0$ is {\em
        independent} (or {\em free}) from $X \cup Y_0$ ($\langle
        X \cup Y_0 \rangle_B$) over the ideal $\Iscr$ in the
        Boolean algebra $B$.

      %  \item If a set $Y$ is {\em independent}  from
%        $X$ over $\emptyset$ then $Y$ is a basis for $\langle X \cup Y\rangle$ over $\langle X\rangle$.

%\item Since $Y$ is independent from $X$ over an ideal $I$ in a
%    Boolean algebra $B$ is expressed by quantifier free formulas,
%    the condition is maintained in any $B'$ extending $B$.
%The next section uses notations and ideas from
\end{enumerate}
\end{observation}
%\sidebar{Do we still want the next?}
From left to right in item 2), %let $\pi$ map $B$ to $B/I$. N
 note that if for any
nontrivial term\footnote{A trivial term (or polynomial) is one which is identically $0$.} $\sigma(\vbar)$, and any $\ybar\in Y$ there is an $a$
with $\sigma(\ybar) \wedge a \not \in I$ then $\pi(\sigma(\ybar)
\wedge a)$ is not $0$ in $B/I$. Conversely, if some
$\sigma(\pi(\ybar)) \neq 0$ then  if $\ybar'$ is  in
$\pi^{-1}(\pi(\ybar))$, then $\sigma(\ybar') \not \in I$.

 %In the construction of complete
%sentence $I$ will be named by a predicate.

\bigskip

The notion of independence above corresponds to the closure system generated by subalgebra \cite[chapter 12]{Gratzer}; it  does {\em not} satisfy the axioms for a matroid (combinatorial geometry); exchange fails.  It is an independence
system (the empty set can be considered independent and subsets of independent sets are
independent.).  But given $X$ and $Y$ independent with $|Y| > |X|$, in general there is no
guarantee that some element of $Y-X$ can be added to $X$ and maintain independence.  But, see Lemma~\ref{freecount}.

The contrast between the notion of independence above and the
following is crucial for the construction here.

\begin{definition}\label{indepset}
Let $X, Y$ be sets of elements from  a Boolean algebra of
    sets. $X$ is independent %(free)
    over $Y$ if for any infinite
    $A$ that is a non-trivial finite  Boolean combination of
    elements of $X$ and any $B$ which is a  non-empty finite
    Boolean combination of elements of $Y$, $A \cap B$ and $A^c
    \cap B$ are each infinite.
\end{definition}

%\sidebar{Is this paragraph right?}
Both kinds of independence will occur in the models in
Section~\ref{zfc+case}. There are models in $\bK_1$, Definition~\ref{defK1}, that are constructed in Construction~\ref{detzfc+}  with a homomorphism from
$P^M_1$ into $\Pscr(P^M_0)$ that does {\em not} transfer from
`independence in the boolean algebra sense' (Definition~\ref{bn}.2 to 'set independence' (Definition~\ref{indepset}.  In
$\bK_2$, there is an isomorphism from $P^M_1$ into $\Pscr(P^M_0)$
that correctly transfers `independence'. (See   Lemma~\ref{nonoise}.)

%\sidebar{Have I reversed the set def and B.A. def - what does each mean?}

\begin{definition}\label{defpushout}
A pushout  consists of an object $P$ along with two morphisms $i_1:
X \rightarrow P$ and $i_2\colon Y \rightarrow  P$ which complete a
commutative square with two given morphisms $f$ and $g$ mapping
an object $Z$ to $X$ and $Y$ respectively such that any morphisms
$j_1,j_2$ from $X$ and $Y$ to a $Q$ must factor through $P$.
\end{definition}

 In \cite{FrGr}, it is
shown by a category theoretic argument that for distributive lattices
the abstract embeddings into the pushout (Notation~\ref{pushout}) are
$1$-$1$ and if $A$ is a Boolean algebra, the images of the embedding
intersect in image of $A$. Thus the variety of Boolean algebras has
{\em disjoint\footnote{Called strong in \cite{FrGr}.} amalgamation}.

We now connect this notion with our version of independence in Definition~\ref{bn}.

\begin{lemma}\label{potoind} Let $D =A \otimes_C B$ be the Boolean algebra
obtained as the pushout (Definition~\ref{defpushout}) of $A$ and $B$
over $C$. Suppose $\Iscr$ is an ideal of $D$ and $I_2\subset A-C$
such that $\langle I_2\rangle_D \cap \Iscr = \emptyset$ and $B -
\Iscr \neq \emptyset$.
 Then $I_2$ is independent from $B$ over $\Iscr$.
\end{lemma}

Proof.  Fix a Boolean-polynomial $p(v_0,\ldots ,v_k)$ (that is not
    identically $0$), and  suppose for contradiction there is a $d \in {B-\Iscr}$ and
     distinct $y_i \in I_2$   with $p(y_0,\ldots ,y_k) \wedge d \in \Iscr$.
    % Now if $p(y_0, \ldots y_k) \wedge a \in \Iscr$, then
     Any morphisms $f_1, f_2$ from $A,B$ to any $D'$ must factor through $D$.  In particular,
     we can extend $\Iscr \cap A$ and $\Iscr \cap B$
to maximal ideals omitting $p(y_0, \ldots y_k)\in A$ and $d \in B$; the resulting map from $D$ that commutes with the induced $f_i$ to the 2-element
     algebra sends all of $\Iscr$ and so $p(y_0,\ldots ,y_k) \wedge d $, but not $d$ or $p(y_0, \ldots y_k)$ to $0$. But there is no such homomorphism.
     $\qed_{\ref{potoind}}$

      % from $D$ which commmutes with both of the $f_i$.S
     %contradicting out choice. %for any choice or $a$, which is impossible.

%
%We can amalgamate Boolean algebras $B$ and $A$ over $C$ by the
%pushout/free product construction of  Notation~\ref{pushout}.
There are several sets of confusing terminology arising from
various perspectives in the study of Boolean algebra and misleading
analogies with, for example, the study of groups. For example,
consider the notion of the product of two Boolean algebras, $A$, $B$.
That is, the structure on the {\em Cartesian  (direct) product} of
$A$ and $B$, obtained by defining the operations coordinate-wise.
Note that, while there are isomorphic copies of $A$ and $B$ in the
product, the natural injections into $A \times \{0\}$, $ \{0\} \times
\{B\} $, map to ideals not sub-Boolean algebras.

A generalization of the dual of the direct product operation is often called the `{\em free product
with amalgamation}'; we will call the free amalgamation of Boolean
algebras $B$ and $A$ over $C$ the one that is  obtained by the
pushout/free product construction of Notation~\ref{pushout}; it is
the coproduct in the category-theoretic language.

%\sidebar{Next almost surely no longer used.}

\begin{notation}\label{pushout} Let $C\subseteq A,B$ be Boolean algebras.  The
disjoint amalgamation $D =A \otimes_C B$ is the Boolean algebra
obtained as the pushout \cite{AvBr} of $A$ and $B$ over $C$. It is
characterized internally by the following condition.  For $a \in A-C,
b\in B-C$: $a \leq b$ in $D$ if and only if there is a $c\in C$ with
$a < c < b$ (and symmetrically). $D$ is generated as a Boolean
algebra by $A \cup B$ where $A$ and $B$ are sub-Boolean algebras of
$D$.
\end{notation}

%Note that this construction changes some atoms to non-atoms: if
%$b,b^- \in B-C$ are finite unions of atoms and $a\in A-C$ is then
%both  of $a \wedge b$ and $a \wedge b^{-}$ must be strictly below $a$
%cannot be $0$ so $a$ is not an atom.

%\sidebar{Is this helpful; we will be amalgamating infinite BA which
%have only finitely many atoms.}
%
%The free amalgam $A \otimes_C B$, where each of $A,B$ have  only
%finitely many atoms must destroy the atomicity of some elements.  (If
%$a$ is atom of $A$ and $b_1, \ldots b_n$ are the atoms of $B$, for at
%least one $i$,  $A \otimes_C B\models 0< a \wedge b_i < a$.) Thus we
%will have to construct a quotient algebra of the free amalgam below
%in order to find an amalgam which does not introduce atoms.

%\sidebar{Remainder of section deleted; amalgamation is done in the
%special situation.}

We will distinguish certain subsets of our models in terms of atoms.
\begin{notation}   An {\em atom} is
an element $a$ of a Boolean algebra such that for every $c$ either $c
\wedge a =a$ or $c \wedge a = 0$.  The element $a$ is a {\em non-trivial atom} if it is neither $0$ nor $1$. For any Boolean algebra $B$,
$\At(B)$ denotes the set of atoms of $B$.

 We  work in a Boolean algebra $P^M_1$ and use $P^M_{4,1}$ for
$\At(P^M_1)$. We will denote by $P^M_4$ the set of {\em finite joins}
of atoms and $P^M_{4,n}$ for those elements that are the join of
exactly $n$ atoms.  $P_4^M$ is always an ideal of $P_1^M$ but it is only a Boolean algebra if it is finite, and even then it will not be a sub-Boolean algebra. A Boolean algebra is {\em atomic}, or in anachronistic terminology, {\em
atomistic} if every element is an {\em arbitrary join} of
atoms\footnote{Equivalently for Boolean algebras, if every non-zero element is above at
least one atom.The conditions are not equivalent on an arbitrary distributive lattice.}.
\end{notation}

For $M$ in the class of finitely generated structures $\bK_0$,  below, the ideal $P^M_4$ will be
atomistic when viewed as a Boolean algebra (with $b^* =1$ and complement as relative compement below $b^*$.) and the
maximal such. For $M$ in the class $\bK_2$ the entire Boolean algebra
$P^M_1$ will be atomistic but this will be false for {\em all} $M$ in
$\bK^{-1}_{<\aleph_0}$ (since it has only finitely many atoms) and for {\em some}  $M$ in $\bK^{1}$ which are not in $\bK^{-1}_{<\aleph_0}$.
We will use the next remark in proving Lemma~\ref{frmodtrans}.

\begin{lemma}\label{batrans} Let $ B_0 \subseteq B_1 \subseteq B_2$ be Boolean
algebras. Suppose $I_i$ for $i<3$  form a sequence of ideals in the
respective $B_i$ with $I_1 \cap B_0 = I_0$ and $I_2 \cap B_1 = I_1$.
If, for $i = 0,1$, $J_i \subset B_{i+1}$ is independent from $B_i$
modulo $I_i$ in $B_{i+1}$, then $J= J_0\cup J_1$ is independent from
$B_0$ modulo the ideal $I_2$.

%\sidebar{sept 21:The text seems to read the last $I_1$ is $I_2$. But I don't
%think that is true.
%sept 23: In fact there were more serious typos.  the independent set should be I 0 and I 1 not I 1 and I 2,
%the final ideal should be I 2. }

\end{lemma}

Proof. Let $\bbar$ be a finite sequence of distinct elements from
$J$. Suppose $\sigma(\ybar)$ is a non-zero term in the same number of
variables as the length of $\bbar$.  For any $d \in B_0 - I_2$, we
must show $\sigma(\bbar)\wedge d \not \in I_2$.
Writing $\sigma$ in disjunctive normal form it suffices to show some
disjunct $\tau$ (which is just a conjunction of literals $y_i$ and
$y_i^-$) satisfies $\tau(\bbar)\wedge d \not \in I_2$. Decompose
$\tau(\bbar)$ as $\tau_0(\bbar_0) \wedge \tau_1(\bbar_1)$ where
$\bbar_i \in J_i$. Since $J_0$ is independent from $B_0$ modulo $I_1$,
$\tau_0(\bbar_0) \wedge d \not \in I_1$ and clearly it is some  $d_1
\in B_1$. Similarly, since $J_1$ is independent from $B_1$ modulo
$I_2$,  $\tau_1(\bbar_1) \wedge d_1 \not \in I_2$. So $ \tau(\bbar)
\wedge d_1 =\tau_0(\bbar_0) \wedge \tau_1(\bbar_1) \wedge d \not \in
I_2$ as required. $\qed_{\ref{batrans}}$

%\sidebar{rest of section added at Kettle Moraine}

\medskip

Although our notion of independence does not satisfy exchange, we are
able to show that under certain conditions each suitable element is a
member of a basis.

% note
%that is a property of sets of elements.  So if $Y$ is independent
%from $X$ over $I$, $B$ is generated as a Boolean algebra by $X \cup Y
%\cup I$ and $Y$ is the disjoint union of $\{Y_i:i<\alpha\}$ then
%letting $B_i = \langle Y_i \cup X \cup I\rangle$ (in $B$),
%\begin{enumerate}
%\item $B = \bigcup_{i<\alpha} B_i$ and
%\item $Y_{i+1}$ is independent of $B_i$ over $I$.
%\end{enumerate}
%
%NOT QUITE -don't want all of $I$.

\begin{lemma}\label{freecount} If $B$ is a countable atomless Boolean algebra, then for any $b \neq 0,1 \in B$, there is a basis $J$ of $B$
that contains $b$.
\end{lemma}

Proof. Observe that by quantifier elimination all non-constant
elements of $B$ realize the same $1$-type.  But then if $A =\langle
a_i\colon i<\omega\rangle$ is a basis for $B$, the automorphism
$\alpha$ of $B$ (guaranteed by $\aleph_0$-categoricity) which takes
$a_1$ to $b$  takes $A$ to $\alpha(A)$ which is a basis containing
$b$. $\qed_{\ref{freecount}}$.

The next result is used in step 2 of the proof of Claim~\ref{cl2}.

%\sidebar{jb I think the next result was implicit but necessary.}
\begin{lemma}\label{basisext} Let $ B_1\subseteq B_2$ be countable Boolean algebras and
suppose $I_2$ is an ideal of $B_2$ and $J_1$ is  a countable  subset
of $B_2$ such that $J_1$ is independent from $B_1$ modulo $I_2$.  If
$b$ is also  independent from $B_1$ modulo $I_2$ and $b \in \langle J_1
\cup I_2\rangle_{B_2}$, then there is a $J_1'$ such that $b \in
J_1'$, $ J_1'$  is independent from $B_1$ modulo $I_2$ and each of $J_1$ and
$J_1'$ generates (with $I_2$) the same subalgebra of $B_2$.
\end{lemma}

Proof.  Let $b^*$ be the image of $b$ when $\pi$ projects $B_2$ onto
$B_2 / I_2$ and $B_3$ denote the image of $\pi(\langle J_1 \cup
I_2\rangle_{B_2})$.  By Lemma~\ref{freecount}, there is a ${J''}_1
\subset B_2 / I_2$ with $b^* \in J^*_1$ that freely generates $B_3$.
Now choose $J_1'$ by choosing a preimage for each element of $J''_1$
and the result follows by Observation~\ref{obs1}.2.
$\qed_{\ref{basisext}}$

\subsection{Defining the Complete Sentence}\label{compsentsec}
\setcounter{theorem}{0}

\setcounter{theorem}{0} In this subsection we construct a complete
$L_{\omega_1,\omega}$-sentence $\phi$, essentially the `existential-completion' of Example~\ref{incompex}. We show in Section~\ref{zfc+case} in an extension of $ZFC$, that $\phi$
has maximal models in $\lambda$ for arbitrarily large $\lambda$ less
than the first measurable cardinal.

% consists of a two-sorted
%structure that consists of  a Boolean algebra on $P_1$ and a
%representation of it as a field of sets on $P_0$ by a predicate
%$R(x,y)$ picking out those $x$ that `belong' to $y$. The basic idea
%is that an extension of a model $M$ by adding an element $a$ to $P_0$
%defines an ultrafilter\footnote{That is, $U$ is an ultrafilter of the
%Boolean algebra with universe $P_1$; it is closed under meet and
%extension. When $R$ is extensional,  the set of $R(P_0,b)$ for $b\in
%U$ is an ultrafilter on $\Pscr(P_0)$. We use both notations.} $U$ on the
%Boolean algebra with domain $P_1$; $U$ is the set of $b$ such that
%$R^{Ma}(a,b)$. If $|M|$ is less than the first measurable $\mu$ this
%ultrafilter cannot be $\aleph_1$-complete.  However, we can construct
%$M$ such that if it has  a proper $(\infty,\omega)$ elementary
%extension, the ultrafilter is $\aleph_1$-complete. For this, we add a
%new predicate $P_2$ and functions $F_n(z)$ so that $P_2$ indexes
%countable families of elements of $P_1$.   %We show here that under additional set theoretic hypotheses
%that any $\aleph_1$-complete ultrafilter on $P_0$ leads in the
%extended vocabulary to a contradiction.
%In the next section we modify the
%construction to obtain the result in ZFC.
%We add a function $G_1$ from $P_0$ into $P_1$ so  that $P_1^M$
%includes an atomic subalgebra $P_4^M$ corresponding to the finite
%Boolean combinations of the singletons of $P_0^M$.
% But taking a
%quotient of $P_1^M$ (by the equivalence relation of `finite symmetric
%difference') is atomless.

Each model
is a member of the class $\bK$ of Example~\ref{incompex}; but Definition~\ref{defK1} describes the  finitely generated models.
This section is devoted to the construction of a countable generic
structure for that class; the details of the construction will be essential for the
main argument in the next section.
Our goal is to build this generic structure as  a Fra{\"i}ss\'{e}-style limit
of finitely generated structures; in each of these structures $P^M_0$
and $P^M_4$ will be finite.

\begin{definition}\label{deftau} $\tau$ is a vocabulary with unary predicates $P_0, P_1,
    P_2, P_4$, binary $R$, $E$, $\wedge, \vee$, %binary function $G_2$,
     unary
    functions $^-$, $G_1$,
   % H_1$, $n$ unary functions $g_{n,i}$ for each $n$,
    constants 0,1 and unary (partial) functions  $F_n$, for $n<
    \omega$.

%
%    \sidebar{july 2: This is probably a dinosaur. $\tau_r$ excludes  the functions $F_n$, for
%    $n\geq r$.}
\end{definition}

We originally introduced the properties of $\bK_{1}$,
% a collection of finitely generated models,
 in two stages ($\bK^{-1}$ and $\bK_{1}$) simply to allow
the reader to absorb the definition more slowly. It turned out in \cite{BaldwinShelahhanfmaxzfc}, that the class $\bK_{-1}$ plays an independent role.

We will study several classes
%\footnote{We originally introduced the properties of $\bK_{1}$,
%% a collection of finitely generated models,
% in two stages ($\bK^{-1}$ and $\bK_{1}$) simply to allow
%the reader to absorb the definition more slowly. It turned out in the later \cite{BaldwinShelahhanfmaxzfc}, that the class $\bK_{-1}$ plays an independent role; so this distinction was wisely made.}
$\bK$ with various subscripts and subscripts. In general for  a class $\bK^i_{<\aleph_0}$  denotes a class of finitely generated structures and either $\hat \bK^i$
or $\bK_i$ denotes the class of all direct limits of models from $\bK^i_{<\aleph_0}$.

We use the word finitely generated in the usual sense. We have a vocabulary with function symbols; each element of a model is given by a term in the finite set of generators. Thus, if $M$ is finitely generated $P^M_0$ and $P^M_2$ must be finite and $P^M_1$ is countable

\begin{definition}\label{defK1}
$\bK^{-1}_{<\aleph_0}$ is the class of finitely generated structures $M$ satisfying.

%0$= \bK^{BA}_0$ if and only if $M$ is a countable
%    $\tau$-structure satisfying condition 1) and is in one of the
%    $K_{0,r}$ in condition 2. $\widehat \bK$ is the
%    class associated with $\bK_0$ as in Definition~\ref{defclass}.
%\begin{enumerate}
%\item $M$ satisfies the following conditions on $P_0$ and $P_1$.
\begin{enumerate}
\item $P^M_0, P^M_1, P^M_2$ partition $M$.
\item $(P^M_1, 0,1, \meet, \vee,< ,^-)$ is a Boolean algebra
    ($^-$ is complement).

%    \sidebar{homomorphism property omitted; was this intentional?}
\item $R \subset P^M_0 \times P^M_1$ with $R(M,b)=
    \{a:R^M(a,b)\}$  and the set of $\{R(M,b): b \in P^M_1\}$ is
    a Boolean algebra. $f^M: P^M_1\mapsto \Pscr(P^M_0)$ by
    $f^M(b) =R(M,b)$ is a Boolean algebra homomorphism into
    $\Pscr(P^M_0)$.

    Note that $f$ is not\footnote{The subsets of $P^M_0$ are {\em
    not} elements of $M$.} in $\tau$; it is simply a convenient
    abbreviation for the relation between the Boolean algebra
    $P_1^M$ and the set algebra on $P_0$ by the map $b \mapsto
    R(M,b)$.

   % \sidebar{Should this be isomorphism or homomorphism? You
%    have the requirement that $R$ is extensional which would
%    makes it an isomorphism.  extensionality in 10/16/15
%    version.}

%\sidebar{Item d) on yellow page 2.2 says $R$ is extensional.
%This would make $f_M$ an isomorphism.  But that requirement
%is violated by the enumeration in case 2 on page 2.14.
%Accordingly, I have only required homomorphism.}
%\sidebar{DT}
    \item $P^M_{4,n} \subseteq P^M_1$ is the set containing each join of $n$
        distinct atoms from $M$; $P^M_{4}$ is the union of the
        $P^M_{4,n}$; $P^M_{4}$ has a maximum element often denoted\footnote{But $b^*$ is not constant in the vocabulary; as the models are extended, $b^*$ changes.} by $b^*$.
         That is, $P^M_4$ is the set of all finite
        joins of atoms (in $P^M_1$).  If $b_1 \neq b_2$ are in $P_4^M$ then $R(M,b_1) \neq
     R(M,b_2)$.

     %   $g^M_n$ is actually a sequence of maps $g^M_{n,i}$ such
%        that if  $P^M_{4,n}(a)$,$\{\langle g^M_{n,i}(a):i<
%        n\rangle$ enumerates the atoms below $a$.

   % $P^M_{4,n} =\{b\in P^M_{1}: |\{c\in P_{4,1}   ^M: c \leq b\}|
%        =n\}$ and $P^M_{4}$ is the union of the $P^M_{4,n}$.

       % \sidebar{$g_n$ omitted; needed? probably-- don't need $\bK_0$ universal.}

      % Let $g_n(\abar) = \bigwedge_{i<n}G^M_1(a_i)$ maps
%        non-repeating sequences  from $(P_0^M)^n$ onto
%        $P^M_{4,n}$, taking all permutations of an $n$-tuple
%        with distinct entries to the same element of
%        $P^M_{4,n}$ and all others to $0$. Do not define
%        $g_n(\abar)$ if the sequence contains repetitions.

   %  \sidebar{$E^M$ omitted since definable but not used.
   %  ????}
   %  \item $E^M$ is an equivalence relation on $P^M_1$ such that
%    $b_1 E^M b_2$ implies $R(M,b_1) \vartriangle R(M,b_2)$ is
%    finite.

\item $G^M_1$ is a bijection from $P^M_0$ onto $P^M_{4,1}$, which by 4)is the non-trivial atoms of $P^M_1$, such
    that $R(M,G^M_1(a)) = \{a\}$.
  %   $H_1^M$ is defined on
%    $P^M_{4,1}$ and is the inverse of $G_1^M$.

   % Note: it is possible $P^M_{4,1}$ and $P^M_{0}$ are empty.

   % \sidebar{SS Mar 15}

\item $P^M_2$ is finite
   (and may be empty). Further, for each $c \in
      P^M_2$ the $F^M_n(c)$ are functions from  $P^M_2$ into $P^M_1$.
      %Note that it is allowed that for all but finitely many $n$, $F^M_n(c) = 0_{P^M_1}$.

%\sidebar{middle of gthe night Sept 27-28:  Put next clause from the Oct 2015 version back in.}

\item If $a \in P^M_{4,1}$ and $c\in P^M_2$ then  for all but finitely many $n$, $ a \nleqslant_M F^M_n(c)$. This implies for each $x \in P^M_0$, $\bigcap_n
\{x\colon(G_1(x)  \in F^M_n(c)\} = \emptyset$.
     \item $P^M_1$ is generated
as a Boolean algebra by $P^M_4 \cup \{F^M_n(c)\mcolon c \in P^M_2, n\in \omega\} \cup X$ where
$X$ is a finite subset of $P^M_1$.

    \end{enumerate}

    %\end{enumerate}
\end{definition}

We denote by $\bK^{-1}$ the class of direct limits of models in $\bK^{-1}_{<\aleph_0}$.
%The $g_n$ were introduced as Skolem functions in  Condition 4
%Definition~\ref{defK1} to guarantee $\bK_{-1}$ was closed under
%%submodel.
%Condition 3 of Definition~\ref{defK1} implies for any $M
%\in \bK_{-1}$,
% any $a \in P^M_0$, $b \in P^M_1$, $M \models R(a,b) \vee
%    R(a,b^-)$.

%\sidebar{beware of redundancy}
%
%The $P^M_1$ of models in $\bK_0$ are  each a free product of a finite
%boolean algebra $B_{n_*}$ and a countable atomless boolean algebra
%$F_\infty$.  $P^M_0$ is in $1$-$1$ correspondence with the atoms of
%$B_{n_*}$, $P^M_2$ indexes families of sequences $F^M_m(c)$ of
%elements of $P^M_1$;
%% if $m<n_*$, $F^M_m(c)\in B_{n_*}$.
% the $F^M_m(c)$ for $m\geq n_*$are independent over the empty set.  The set $\{F_m(c): m
%\geq n_*, c \in P^M_2\}$ are a basis  for $F_\infty$. The $F^M_m(c)$
%for $m< n_*$ are Boolean combinations of elements in $B_{n_*}$ and
%the $F^M_n(c)$ with index greater than $n_*$;  as the models are
%extended and $n_*$ grows all the $F^M_n(c)$ eventually have atoms
%below them. (See Lemma~\ref{nonoise}.)
%
 We now add requirements to Definition~\ref{defK1}, to ensure that no elements of $P^M_1$ are needed as generators and to lay the ground for the study of free extensions. (See Definition~\ref{deffreesub}.)
We  refine the class $\bK^{-1}_{<\aleph_0}$ from  Definition~\ref{defK1} to a
 class $\bK^{1}_{<\aleph_0}$; here  the structure is witnessed by a family of
  witnesses  $\langle n_*, \bB, b^*\rangle$. The class of direct limits of these finitely generated structures generate will be denoted
  $\bK^{1}$. From $\bK^{1}$ we will derive the rich class
 $\bK_2
=\Krich$ in Definition~\ref{richname}.

\begin{definition}\label{k0}
  $M$ is in the class of structures $\bK^{1}_{<\aleph_0}$ if  $M \in \bK^{-1}_{<\aleph_0}$ and there is a {\em witness}
  $\langle n_*, \bB, b^*\rangle$
  such that:

  \begin{enumerate}
  \item $b^* \in P^M_1$ is the supremum of the finite joins of
      atoms in $P^M_1$. Further, for some $k$,\\ $\bigcup_{j\leq k}
      P^M_{4,j} = \{c: c \leq b^*\}$ and for all $n>k$,
      $P^M_{4,n} =\emptyset$.
  \item $\bB = \langle B_n: n \geq n_*\rangle$ is an increasing
      sequence of finite Boolean subalgebras of $P^M_1$.

 % \sidebar{DT   was $B_{n_*} = P^M_4$:}

\item   $B_{n_*}\supsetneq    \{b\in P^M_1: b \leq b^*\} =
    P^M_4$; it is generated by the subset $P^M_4
    \cup\{F^{M}_n(c)\colon n< n_*, c \in P^{M}_2\} $.

    %\sidebar{following added by SS in Mar 2017. Is it necessary?}

Moreover, the Boolean algebra $B_{n_*}$ is free over the ideal
$P^M_4$  (equivalently, $B_{n_*}/P^M_4$ is a free Boolean
algebra\footnote{A further equivalence:
 $|Atom(B_{n_*})| -
|P^M_{4,1}|$ is a power of two.}.)

%\item If $n\geq n_*$, then $B_n$ is generated as a Boolean algebra by  $B_{n_*}\cup \{F^M_i(c) : c\in P^M_2, n_* \leq i\leq n $.

  \item $\bigcup_{n\geq n_*}B_n = P^M_1$.
  \item %$P^M_2$ is finite and not empty. Further,
  For each $c \in
      P^M_2$ the $F^M_n(c)$ for $n< \omega$ are distinct and independent over
       $\{0\}$.

      %; the $F^M_m(c)$ for $m < n_*$ are in $B_n$.

  %\sidebar{previous seems included in next}

% \item   If $c_0,\ldots c_{k-1} \in P^M_2$ are pairwise
%    distinct then
%$$\{F_{n}(c_\ell): \ell < k, n\geq n_*\}$$
%has no repetitions and for $n\geq n_*$,

% \sidebar{RESTORED SEPT 24 AM SCAN page 11/36 printed page 5:deleted in Sept 21 scan:$F_n(c_i) \meet b^* = 0$ (in
% the Boolean algebra on $P^M_1$).} \label{fixR}

%\sidebar{DT   was $B_{n_*} = P^M_4$: mar 13 - rest not right yet}
\item  The set $\{F_m(c):m\geq n_*, c \in P^M_2\}$ (the
    enumeration is without repetition) is free from $B_{n_*}$
    over\footnote{As in Definition~\ref{bn}.2 with $X=
    \emptyset$.}  $\{0\}$. $B_{n_*} \supsetneq P^M_4$ and $F_m(c)
    \wedge b^* = 0$ for $m\geq n_*$. (In this definition, $0 =
    0^{P^M_1}$.)

%\sidebar{DT  tried to fix the abysmal $x_c$ notation.}
 In detail,
let $\sigma(\ldots x_{c_i} \ldots)$ be a Boolean algebra term in
the variables $x_{c_i}$ (where the $c_i$ are in $P^M_2$)  which is
not identically $0$. Then, for finitely many $n_i\geq n_*$ and a finite sequence of  $c_i\in P^M_2$:
$$ \sigma(\ldots F_{n_i}(c_i) \ldots) > 0$$
and some $n<\omega$. Further, for any non-zero $d \in B_{n_*}$ %P^M_1$
with $d \wedge b^* =0$,  (i.e. $d \in B_n - P^4_M$),
$$
\sigma(\ldots F_n(c) \ldots) \wedge d > 0.$$

%\sidebar{Where used?}

\item For every $n\geq n_*$, $B_n$,
% and so $P^M_1$,
%\sidebar{$B_n$ was written.}
is generated by $B_{n_*} \cup \{F_m(c):n > m\geq n_*, c \in
P^M_2\}$. Thus $P^M_1$ and so $M$ is generated by $B_{n_*}\cup
P^M_2$.

%\item If $n< n_*$ and $c \in P^M_2$, $F^M_n(c) \in B_{n_*}$
%\item $P^M_4 \subseteq B_{n_*}$.
\end{enumerate}

%and $\langle F_n^M(c_i): c_i \in P_2^M: p<i
%\rangle$ is independent over $\langle F_n^M(c_i): c_i \in P_2^M:
%i<p\rangle\cup P_4^M \cup X$.

\end{definition}

%Condition 7 implies  for $n\geq n_*$ that the $F_n(c)$ are independent over $P^M_4$.
%Conditions 1 and 3 imply $P^M_4 \subseteq B_{n_*}$.
%In fact, $B_{n_*}$ could be taken to be generated by $P^M_4 \cup \{F_m(c): c< n_*, c \in P^M_2\}$.

\begin{remark}\label{neededfacts}

{\rm
The first part of Condition 6 of Definition~\ref{k0} implies
% that if $a\in P^M_0$ and
%$c\in P^M_2$ then for every large enough $n$, $a \not \in
%R(M,F_n(c))$. That is,
condition 8 of Definition~\ref{defK1}. The second part of condition 6
implies, in particular, that if $b \in P_1^M-P_4^M$, there are
infinitely many elements below $b$ in $P_1^M$.  Note that the free generation condition of 6) is not preserved by arbitrary direct limits; in particular it will
fail in the $P_0$-maximal model of cardinality $\lambda$. However, our construction in an extension of ZFC of $\bK_1$-free extensions   (Definition~\ref{deffreesub}) will guarantee the  $\bK_1$-freeness of submodels that are less than $\lambda$-generated.
% is met
%in a very strong way. Equivalently $\bigcap_n R(M,F^M_n(c)) = \emptyset$.

Note that if $\langle n_*, \bB, b^*\rangle$ witnesses $M \in \bK^1_{<\aleph_0}$
then for any $m\geq n_*$, so does $\langle m, \bB, b^*\rangle$.}
\end{remark}

The following lemma shows the prototypical models in $\bK^1_{<\aleph_0}$ in fact
exhaust the class. Note that each $P^M_1$ is an atomic Boolean algebra.

\begin{lemma}\label{decomp} For any $M\in \bK^1_{<\aleph_0}$, $P^M_1$ has a natural decomposition as a
product of an atomic and an atomless Boolean algebra. %\footnote{The
%atomic component need not have a maximal element (in $\Krich$).}.
%% and certainly $1\in P^M_1$ is not
%in $P^M_4$.
\end{lemma}

Proof. Let $M\in \bK^1_{<\aleph_0}$, witnessed by  $\langle n_*, \bB,
b^*\rangle$. Then the atomic part, $P^M_4$, is the collection of
elements of $P^M_1$ that are $\leq b^*$. And the independent
generation by the $F_n^M(c_i)$ for $n\geq n_*$ and $ c_i \in P_2^M$
shows the quotient $P^M_1/P^M_4$ is atomless. $\qed_{\ref{decomp}}$

Condition Definition~\ref{k0}.3 guarantees:

\begin{lemma}\label{fgk0} Each structure in $\bK^1_{<\aleph_0}$ is finitely generated  by $P^M_0  \cup P^M_2$.
\end{lemma}

%Proof. Let $M\in \bK^1_{<\aleph_0}$, witnessed by  $\langle n_*, \bB,
%b^*\rangle$. Then $M$ is generated by $P^M_0  \cup P^M_2$.
%$\qed_{\ref{fgk0}}$

\medskip

\begin{lemma}\label{fewmodels} $\bK^1_{<\aleph_0}$ is countable.
\end{lemma}

Proof. Let $M$ be in $\bK^1_{<\aleph_0}$, witnessed by  $\langle n_*, \bB,
b^*\rangle$.  The isomorphism type of $M$ is determined by the
structure on $P^M_4$ induced by the  $F_n(c_i)$
 %for $n<
%n_*$
and $c_i \in P^{M}_2$.  If $m \geq n_*$, $F_m(c_i) \wedge b^* = 0$ so
they leave no trace on $P^M_4$. Since this tail, $\{F_m(c_i) \colon m
\geq n_*\}$ generates an atomless boolean algebra in the sense of
$P^M_1$, that boolean algebra is $\aleph_0$ categorical. But there can be only
countably many structures induced on the finite $P^M_4$ by the
countable set $F_n(c_i)$ through the formulas $x <F_n(c_i)$ which
determine the values of $R$ on $P^M_4$ since only the $F_m(c_i)$ for
$m < n_*$ have non-empty intersection with $P^M_4$ (i.e. are above atoms) and $P^M_2$ is finite.
$\qed_{\ref{fewmodels}}$

%The  sEpt visit  idea of using declaring all the $F_N(c)$ do not intersect $P^{n_*}_4$ was abandoned on Sept 21.
%But condition~\ref{fixR} guarantees that for each $i$, $F_n(c_i) \meet b^* = 0$,
%so for $n\geq n_*$, $R(M,F_n(c_i))$ is empty.

\medskip

%\sidebar{It appears Lemma~\ref{decomp} isn't used.}

%This decomposition is not claimed to be either {\em uniformly }
%definable or unique.  It is not given by `symmetric difference is
%finite' on the finitely generated models.
%  That equivalence relation
%has only one class on $\bK_0$ but many on models in $\Krich$.
%
%\sidebar{CHECK LAST PARAGRAPH}

\begin{definition}\label{k1def}  The class $\bK_1= \hat \bK^1_{<\aleph_0}$ is the collection of all
direct limits of models in $\bK^1_{<\aleph_0}$.
\end{definition}

%Recall that members of $\bK^1_{<\aleph_0}$ must have $P_2 \neq \emptyset$ but $P_0$ may be empty.

\begin{lemma}\label{minmod} There is a minimal model $M_{min}$ of $\bK^1_{<\aleph_0}$,
that can be embedded in any model of $\bK_1$.
\end{lemma}

Proof. Let $P^{M_{min}}_0$  be empty,  so $P^{M_{min}}_4 =\{0\}$. Also, let
$P^{M_{min}}_2$ be empty.
%
%contain a single element $c$,  and let the $F^{M_{min}}_n(c)$
%be independent generators of $P^{M_{min}}_1$. Since $M\in \bK^1_{<\aleph_0}$ implies
%$P^{M}_2 \neq \emptyset$, for the embedding we can map $c$ to any element $c_0$ of $P^{M}_2$
%and $F^{M_{min}}_n(c)$ to $F^{M}_n(c_0)$ by Definition~\ref{k0}.5.

%\sidebar{This is an abandoned notion of minimal.}
% We say $M$ is free
%over the emptyset or
%    simply free if $M$ is free over $M_{min}$.; let the $F^{M_0}_n(c)$ be independent
%generators of $P^{M_0}_1-P^M_4$; $P^M_4$ is empty.
%$\qed_{\ref{minmod}}$
%
%
%\medskip Thus, in the minimal model $P^{M_0}_1$ is the direct product
%of the 2-element Boolean algebra with an atomless Boolean algebra.
%More generally, if $M \in \bK_0$, $P^{M}_1$ is the direct product of
%a Boolean algebra whose points $<1$ form the ideal $P^M_4$ with an
%atomless Boolean algebra.}
%\sidebar{DT}

%\sidebar{Since there were a couple of quite misleading typos in Lemma
%3.2.10 (Oct 2016), I include the corrected version with proof.}

\begin{lemma}\label{fixbase}
If $M_0 \subseteq M_1$ are both in $\bK^1_{<\aleph_0}$, witnessed by $\langle
n^i_*, \bB^i, b^i_*\rangle$, for $i=0,1$, then for sufficiently large $n$, $B^0_n
= B^1_n \cap P^{M_0}_1$.
\end{lemma}

Proof. Recall $\bB^i= \langle B^i_n\colon n< \omega\rangle$. Since
the $B^1_n$ exhaust $P^{M_1}_1$, $B^0_{n_*}$ is finite, and for $c
\in P^{M_0}_2$ and all $r$, $F^{M_1}_r(c) =F^{M_0}_r(c)$, for all
sufficiently large $n$, $B^1_n$ contains the $F^{M_0}_r(c)$ for $r<n$
and thus $B^0_n$. But if some $b\in B^1_n \cap P^{M_0}_1 $, but  is
not in  $B^0_n$ then for some $k$, $b\in B^0_{k+1}-B^0_{k}$. But then
$B^0_{k+1}$ is not generated by $B^0_{n^*}$ along with the
$F^{M_0}_r(c)$ for $r < k$. $\qed_{\ref{fixbase}}$

Note that if the conclusion of Lemma~\ref{fixbase} holds for $n$, it holds for all $m\geq n$.

\medskip

We now introduce some  special notation for this paper by
defining {\em $\bK_1$-free over} ({\em $\bK_1$-free extension} of) for models in
$\bK_1$.  $M_2$ is a $\bK_1$-free extension of $M_1$ if not only is
the image of $P^{M_2}_1$ in the Boolean algebra $P^{M_2}_1/P^{M_2}_4$
 a free extension of the image of $P^{M_1}_1$ but the $F_n(c)$ satisfy
technical conditions which allow the preservation of this condition under unions of chains.

\begin{definition}\label{deffreesub} When  $M_1 \subseteq M_2$ are both in
 $\bK_1$, we say $M_2$ is {\em $\bK_1$-free over} or is a {\em $\bK_1$-free extension of } $M_1$ and write $M_1
        \subseteq_{fr} M_2$, witnessed by $(I,H)$ when
%\begin{enumerate}
%\item  We say $M_1 \subseteq_* M_2$ if
%\begin{enumerate} \item $M_1 \subseteq M_2$ are both from $\bK_1$. %$\widehat{\bK}$.
%%\item If $c \in P_2^{M_2} \setminus P_2^{M_1}$, then for some
%%    $r$, $\langle F^{M_2}_n(c)\colon c\geq r\rangle$ is
%%    independent in $P_1^{M_2}$ over  $P_1^{M_1}\cup P_4^{M_2}
%%    $.
%%    \end{enumerate}
%   \item The freeness is witnessed by $(I,H)$ when

       \begin{enumerate}
%\item   $M_1 \subseteq_* M_2$ are from $\widehat{\bK}$.
 \item $I \subset   P^{M_2}_1 - (P^{M_1}_1 \cup P^{M_2}_4)$ satisfies i) $I \cup P^{M_1}_1 \cup P^{M_2}_4$ generates
     $P^{M_2}_1$ and ii) $I$ is independent from
$P^{M_1}_1$ modulo $P^{M_2}_4$ in $P^{M_2}_1$. (Definition~\ref{bn}.2.)

% \item There is a function $H$ from $P^{M_2}_2$ into $\omega$
%     such that
%
%   % is a free basis for $P_1^{M_2}$ over  $P_1^{M_1} \cup
%    P_4^{M_2}$ (equivalently over $P_1^{M_1} \cup P_0^{M_2}$).
   \item There is a function $H$ from $P_2^{M_2} \setminus
       P_2^{M_1}$ to $\NN$ such that the $F_n(c)$ for $n\geq
       H(c)$ are distinct and
        $$ \{F_n^M(c)\colon  c \in P_2^{M_2} \setminus P_2^{M_1}
     \text{ and } n \geq H(c)\} \subset I$$
     and for every $c\neq d \in P^{M_2}_2$, $\{n\colon (\exists m) F^{M_2}_m(c) =  F^{M_2}_n(d)\}$ is finite.

 %    index distinct elements of  $I$.

%\sidebar{H is needed because in order to realize a specific sets
%$C$ as $R(M,b)$, need to have models that are not nearly free.}
%
%
%        satisfying
%      %  \begin{itemize}
%         if $c_1 \neq c_2$ are from  $P_2^{M_2} \setminus
%            P_2^{M_1}$ and $n_1 \geq H(c_1)$, $n_2 \geq
%            H(c_2)$ then $F_{n_1}^{M_2}(c_1) \neq
%            F_{n_2}^{M_2}(c_2)$

               %\end{itemize}
\end{enumerate} %\end{enumerate}
 % \item
  We say $M$ is
 $\bK_1$-free over the empty set or
    simply $\bK_1$-free
    if $M$ is a $\bK_1$-free extension of  $M_{min}$.

%   $M$ is free if it is free over the empty model i.e., $P^M_1$ has a free basis over $P^{M}_4$.
%                \item $M$ is weakly free over $N$ if $M$ contains a free extension $M'$ of $N$.
               % \end{enumerate} \end{enumerate}

    \end{definition}

%We also say $M$ is {\em nearly free} if $P^M_1/P^M_4$ is a free Boolean algebra on infinitely many generators and  $M$ is {\em nearly atomless} if $P^M_1/P^M_4$ is an atomless Boolean algebra. These notion are synomomous only if $P^M_1P^M_4$  is countably infinite.

    %Note that if $M\in \bK_{0,r}$, $M$ is free since $P^M_4$ is generated by
%    $P^M_0$. Moreover, if $N_1$ is  $(M,k)$-closed then there exists an
%    $\check M$ (generated by the $F_n(c)$ with $c\in N_1$ and $n<k$ such that
%$N_1 \subseteq_{fr} \check M \subseteq M$.

%\begin{remark}  Note that a free model in $\bK_1$, like those in
%$\bK_0$, has %is atomistic,there is
%an ideal picked out by $P^M_4$, and $P^M_1/P^M_4$ is an atomless
%boolean algebra.
%\end{remark}
%
%We need only speak  of the generation of $P^M_1$ and  $P^M_2$  since
%for  $M \in \bK_1$  elements of $P^M_0$  are generated from $P_1$ by
%$G_1$. Here are some  basic facts about free extensions.

%\sidebar{Shelah sept 21 page 1.24 has a couple of more consequences;
%are they used?}
\begin{lemma}\label{frmodtrans}
\begin{enumerate}
\item If $M_1 \subseteq_{fr} M_2$ by $(I_1,H_1)$ and $M_2
    \subseteq_{fr} M_3$ by $(I_2,H_2)$ then $M_1 \subseteq_{fr}
    M_3$ by $(I_1 \cup I_2, H_1 \cup H_2)$. Thus,
    $\subseteq_{fr}$ is a partial order.
\item More generally, if $M_\alpha$ with $\alpha < \delta$ is
    a continuous $\subseteq_{fr}$-increasing sequence       then $M= \bigcup
    M_\alpha$ satisfies $M_\alpha \subseteq_{fr} M$ witnessed by
    $(\bigcup_{\alpha < \beta < \delta}I_\beta, \bigcup_{\alpha < \beta < \delta}H_\beta)$.
\end{enumerate}
\end{lemma}

Proof. By Lemma~\ref{batrans} (taking the ideals as $P^{M_2}_4$ and
$P^{M_3}_4$), $I_1 \cup I_2$ is free from $P^{M_1}_1$ over
$P^{M_3}_4$.
% But we can take $H_3(c)$ as $H_i(c)$ for $i = 1,2$ (note
$H_1 \cup H_2$ is well-defined since the $H_i$ are defined on disjoint sets.
    %because the $F^{M_{i+1}}_n(c)$ for $c \in P^{M_{i+1}}_2 - P^{M_{i}}_2$
    %cannot lie in $P^{M_{i}}_1$ by the independence.
Part 2 follows by induction. Successors are similar, while limits
 are automatic.
 $\qed_{\ref{frmodtrans}}$

\begin{remark}\label{monnote} In an increasing chain such as that of
Lemma~\ref{frmodtrans}.2, if some $b \in P^{M_{\alpha+1}}_1$ is free from
$P^{M_{\alpha}}_1$ modulo $P^{M_{\alpha+1}}_4$  then $b$ is also free
from $P^{M_{\alpha}}_1$ over $P^{M_\beta}_4$ for any $\beta > \alpha$
since $P^{M_\beta}_4 \cap P^{M_{\alpha+1}}_1 = P^{M_{\alpha+1}}_4$.
\end{remark}
% the requirements on the $I_i$ follow directly
%Lemma~\ref{batrans}  (taking the ideals as $P^{M_2}_4$ and
%$P^{M_3}_4$) and the given generating conditions.   And clearly the
%$F_n(c)$ are contained in $I_1 \cup I_2$ and continue to be
%independent.

%Note that each model in $\bK_0$ is free over the emptyset.
%\sidebar{CHECK NEXT CAREFULLY and add paragraph on why it is the MAIN
%point. revised Claim 1.13 on page 1.26 of Sh sept 21 (previous
%version used closed). The corollary is called straightforward on page
%1.28 of sept 21. The proof there says only asks to verify a
%triviality.}

The next lemma uses the requirement that the $B_n$ in the witnessing
sequence are free Boolean algebras.

\begin{lemma}\label{subeqfree} If $M_0 \subset M_1$ are both in $\bK^1_{<\aleph_0}$ then $M_0 \subset_{fr}
M_1$.
\end{lemma}

Proof. %Choose  $n_*$ as the maximum of $n^i_*$ for $i<2$;
We can
assume by Lemma~\ref{fixbase} that the $n^i_*$ for $i<2$ are equal
and that $B^1_{n^1_*} \cap P^{M_0}_{1} = B^0_{n^0_*}$.  Since the
$B_{i,n^*}$ are free from $P_4^{M_i}$ over $\emptyset$, we can choose bases $I_0, I_1$
for $B_{0,n^0_*}$ and $B_{1,n^1_*}$ respectively. Now $I_{0} \cup I_1 \cup
\{F^{M_i}_n(c)\colon i<2, n \geq n^i_{n^i_*}, c \in P^{M_i}_2\}$ is a
free basis of $P^{M_i}_1$ over $P^{M_i}_4$. Hence $(I_2 \setminus
I_1) \cup\{F^{M_i}_n(c)\colon i<2, n \geq {n^i_*}, c \in
P^{M_i}_2\}$ is the required $I$ from Definition~\ref{deffreesub}
with $H(c) = n^*$ for all $c$. $\qed_{\ref{subeqfree}}$

\medskip
Crucially, Lemma~\ref{subeqfree} fails in general if $\bK^1_{<\aleph_0}$ is replaced by
$\bK_1$.  Lemma~\ref{subeqfree} immediately yields.

\begin{corollary}\label{kotofree} Each model $N$ in $\bK^1_{<\aleph_0}$ is $\bK_1$-free over the empty set.
\end{corollary}

%Proof. Fix $c_0 \in P^N_2$ as the generator of an embedding of the
%minimal model. Let $\langle n_*, \bB, b^*\rangle$ witness that $N \in
%\bK_0$. Then $B_{n_*}$ has a free basis $I_1$ over $P^{N}_4$ by
%Definition~\ref{k0}.3. Moreover,  by Definition~\ref{k0}.6,  $(I,H)$
%where $I = I_1 \cup \{F_n(c)\colon c \in P^{N}_2 , n \geq n_*\}$ and
%for each $c\in P^N_2$ setting $H(c) = n_*$ witnesses that $M$ is a
%free extension of $M_{min}$. $\qed_{\ref{kotofree}}$

To find large $\bK_1$-free models we apply
Lemma~\ref{frmodtrans}.2  to construct a sequence of $\bK_1$-free extensions. We now show that if $M_1$ is $\bK_1$-free, $N_1 \subseteq M_1$ and $N_1 \subseteq N_2$ with $N_2$ a finitely generated extension  of the  finitely
generated substructure $N_1$, then $M_1$ and $N_2$ can be amalgamated over $N_1$.  Note that by Lemma~\ref{subeqfree}, on $\bK^1_{<\aleph_0}$, $\subseteq$ is the same as $\subseteq_{fr}$. There are three key ingredients in the amalgamation proof: $N_1$ and
$N_2$ must be finitely generated; this is reflected positively in the ability to
employ the witnessing sequences $\bB^i$ in the proof but also by  the key role
in the proof of the finite set $P^{N_2}_{4,1}- P^{N_1}_{4,1}$. Secondly,    $M_1$
must be $\bK_1$-free.
%This is used in the proof that $M_2 \in \bK_1$.
 Thirdly, we must ensure that `atomicity' is preserved in constructing extensions of Boolean algebra so the definitions of $P_4$ and $P_{4,1}$ are `absolute' between models.  It is this third condition which drives the complexity of steps 1 to 3 in the following proof. The free amalgam $D= A \otimes_C B$, where either of $A,B$
has only finitely many atoms must destroy the atomicity of some
elements. (If $a$ is an atom of $A$ and $b_1, \ldots b_n$ are the
atoms of $B$, for at least one $i$, $A \otimes_C B\models 0< a \wedge
b_i < a$.) Thus we will have to construct a quotient algebra of the
free amalgam in step 3 below in order to find an amalgam which does not
destroy atoms.

%We have
%several simple methods.

%\sidebar{jb Apr 10 -- commented out my earlier `proof' of 1.19}
%\sidebar{JB Apr. 7. Here is a proposed proof of 1.19 in the recent
%scans.  It is just an observation that free extension and finitely
%generated extension was the essence of the amalgamation lemma.
%
%It that is actually correct, I could state a somewhat more general
%theorem.  I worry about your proof of 1.19 that it does not seem to
%take account of preserving the atoms of $M_2$ in the notation below.
%}
%
%Now we note that the rather than `both models are finitely generated'
%the key ingredients in the amalgamation result are the three models
%are all free and the intersection model is finitely generated.

%\newpage
% $P^M_4 \cup \{\langle
%1,1\rangle, \langle a,1\rangle, \langle b,1\rangle, \langle
%0,1\rangle,\}$.

%
%
%\sidebar{This argument shows there is a problem when a factor is
%finite.  We don't have that any more unless by accident
%$P^{M_1}_{4,1}$ add only finitely many points -- so is this still
%relevant?
%
% The free amalgam $D= A \otimes_C B$, where either of $A,B$
%has only finitely many atoms must destroy the atomicity of some
%elements. (If $a$ is ab atom of $A$ and $b_1, \ldots b_n$ are the
%atoms of $B$, for at least one $i$, $A \otimes_C B\models 0< a \wedge
%b_i < a$.) Thus we will have to construct a quotient algebra of the
%free amalgam below in order to find an amalgam which does not
%introduce atoms.}

%\sidebar{CAn free be nearly free or is this stronger?}

 \begin{theorem}\label{amalfree}    Suppose $M_1\in \widehat \bK =\bK_1$ is $\bK_1$-free and $N_1 \subset M_1$.
 Let $N_1 \subset N_2$ with both in $\bK^1_{<\aleph_0}$.

Then there are an $M_2$ and an $f$ such that:

    \begin{enumerate}
    \item $M_2 \in \bK_1$, $M_1 \subseteq_{fr} M_2$ and so $M_2$
        is $\bK_1$-free.  %Similarly $N_2 \subseteq_{fr} M_2$.
    \item $f$ maps $N_2$ into $M_2$ over $N_1$. Moreover, the
        image in $M_2$ of $N_2$ is $\bK_1$-free over $N_1$.

    \end{enumerate}

   \end{theorem}

   %\sidebar{stopped in middle apr 10 afternoon: apr 13, trying to copy Saharon's proof.}

Proof.  We lay out the situation in more detail. $M_1$ is $\bK_1$-free means
that $M_1$ is $\bK_1$-free over $M_{min}$ by $(I_1,H_1)$. For $i = 1,2$, let
$\langle n^i_*, \bB^i, b^i_*\rangle$ witness that $N_i \in \bK^1_{<\aleph_0}$.
Suppose $N_1 \subseteq_{fr} N_2$ is witnessed by $(I_2,H_2)$.
Invoking Lemmas~\ref{fixbase} and \ref{neededfacts}, we can rename
$n^i_*$ and rechoose $n_*$ for $N_2$ so  that $n^1_* = n^2_* =n_*$ and $B^1_n =
B^2_n \cap N_1$ for $n \geq n_*$, and (since $P^{N_2}_2$ is finite)
for each $c \in P^{N_1}_2$, $H_1(c) \leq n_*$. Let $J_1 \subset
B^1_{n_*}$ be the pre-image of the basis of $B^1_{n_*}/P^{N_1}_4$.
Then, since $J_1/P^{N_1}_4$ is a generating set of
$B^1_{n_*}/P^{N_1}_4$, for each $b \in B^1_{n_*}$, there is a Boolean
combination  $b'$ of elements of $J_1$ such that $b' \triangle b \in
P^{N_1}_4$.  Note also, that by our choice of $n_*$
(Definition~\ref{k0}.6), if $b \in P^{N_1}_1$ is above an atom of
$P^{N_2}_1$, $b \in B^1_{n_*}$.
Let $k = |P^{N_2}_{4,1}-
P^{N_1}_{4,1}|$, fix $a_0 \ldots a_{k-1}$ listing a new set $A$, and
let {\em $f$ be $1$-$1$ function from $P^{N_2}_{4,1}- P^{N_1}_{4,1}$ onto
$A$}; $A$ contains an image of each new atom in $N_2$.
%\sidebar{This is a different concept than the dimension of
%$B_{2,n_*}$ over $P^{N_2}_4$; I don't know how they are related or if
%it matters. }

%Mistaken attempt to prove only finitely generated amalgamation-- back to mar 2020 main idea of proof.
%Let $\mathbb{B}_1$ be the Boolean algebra which is the  free product with amalgamation: $P^{N_1}_1 *_{P^{N_0}_1}P^{N_2}_1$ and $\ell<3$, let $h^1_
%i$ be the natural embedding of  $P^{N_i}_1$ into  $\mathbb{B}_1$.
%
%This construction may destroy atoms. If $a,b$ are atoms  in $N_1- N_0$ and $N_2- N_0$, it may well be that  $\mathbb{B}_1 \models 0< h^1_1(a) \meet h^1_2(b) < h^1_1(a)$. So we shall divide $\mathbb{B}_1$ by an ideal.
%
%We may assume the ${P^{N_i}_1}$ are disjoint. We write $\b2$ for the 2-element Boolean algebra. For each atom $a \in A = P^{N_1}_{4,1} \cup P^{N_2}_{4,1}$, and each $i< 3$ we define a map $g^i_a: P^{N_i}_1 \rightarrow \b2$ such that:
%{$\bigoplus ^i_a$} \nolinebreak
%\begin{enumerate}
%\item $g^i_a$ is a homomorphism from ${P^{N_i}_1}$ to $\b2$;
%\item $g^0_a \subseteq g^i_a$ for $i=1,2$;
%\item for each $c\in {P^{N_i}_2}$ here are only finitely many $n$  with $g^i_a(F^{N_i}_n) =0$.
%    \end{enumerate}
% For this, if $a\in P^{N_i}_{4,1}$ with $i=1,2$ and $b\in P^{N_i}_1$ let
%
%
%   \begin{equation}
%    I(T,\aleph_\alpha)
%    \begin{cases}
%     = 1 \, & a \leq_{N_i} b \\
%    0 &  a \meet_{N_i} b = 0
%    \end{cases}
%  \end{equation}
%
%  Letting $g^0_a = g^i_a$, i) and ii) are clear.  $\bigoplus.3 holds by
%Lemma~\ref{k0}.5. We extend each $g^i_a$ to $\mathbb{B}_1$
%
%*********************

{\bf Step 1:} {\em Construct
%a Boolean algebra $\mathbb{B}_1$ that has exactly the required atoms.}  More precisely,
%construct
 a Boolean algebra $\mathbb{B}_1$
%that \sidebar{includes// SS writes is generated by}
that
 is generated by
$P^{M_1}_1 \cup A$
 and so that the atoms of $\mathbb{B}_1$ are   $P^{M_1}_{4,1} \cup A$.}
  For this demand,  let $\Dscr_\ell$, for each $\ell < k$, be an ultrafilter of the Boolean algebra $P^{M_1}_1$,
disjoint from $I_1-J_1$ such that for $b \in P^{N_1}_1$,  $b \in
\Dscr_\ell$ if and only if $N_2 \models f^{-1}(a_\ell) \leq b$. (Such
an ultrafilter exists as the set $\{b \in P^{N_1}_1\colon
f^{-1}(a_\ell) \leq b\}$, as noted in last paragraph, contains no element of $I_1-J_1$ and  is a filter on
  $P^{N_1}_1$ that can be extended to an ultrafilter on the Boolean algebra $P^{M_1}_1$.)     % By
%a similar argument, there is no $b \in P^{M_1}_{4}$ with $N_2\models
%a_i \leq b$.

    Now let $X$ be the union of the Stone space of $P^{M_1}_1$, denoted $S(P^{M_1}_1)$ with $A$.
% with new
%elements $a'_i$ for $a \in A$
For $b \in P^{M_1}_1$, let $$X_b =\{d \in S(P^{M_1}_1): b \in d\}
\cup \{a_\ell: b \in \Dscr_\ell\}.$$

Now let $\mathbb{B}_1$ be the subalgebra of $\Pscr(X)$ generated by
the $\{X_b: b \in P^{M_1}_1\} \cup A$.  Now, generalizing the Stone
representation theorem, we embed  $P^{M_1}_1 \cup A$ into $\bbB_1$ by a map $g$; let $g(b) =X_b$  for $b\in P^{M_1}_1$
and $g(a) = \{a\}$ for $a\in A$.

%\sidebar{??? Abusing notation, we write
%$f(a)$ for $\{f(a)\}$.}

Since $P^{N_2}_{4,k} \cap P^{N_1}_{1}  =P^{N_1}_{4,k}$,
there can be no non-zero $b \in P^{N_1}_{4}$ and so no non-zero $b \in P^{N_2}_{4}$  with $N_2\models  b < f^{-1}(a_i)$.
%So the elements of $A$ are atoms in $\mathbb{B}_1$.
  Note  i) that  for  $b\in P^{M_1}_1$, $b \in
\Dscr_\ell $ iff $\mathbb{B}_1 \models f^{-1}(a_i) \leq b$   and ii)  that $e$ is an atom of $P^{M_1}_1$ if and only if $X_e$ is a principal ultrafilter in $\mathbb{B}_1$.
Thus, the atoms of $\mathbb{B}_1$ are exactly
$P^{M_1}_{4,1} \cup A$.

%We imbed $P^M_1$ into $\mathbb{B}_1$ by identifying $b$ with
%$X_b$ and identify $P^M_1$ with its natural image.
%
{\bf Step 2:} {\em Find a  sub-Boolean algebra $\bbB^*$ of $\mathbb{B}_1$ that is
a suitable  base for amalgamating $\bbB_1$ with $P^{N_2}_1$}. For this, denote by $\mathbb{B}^*$ the
sub-Boolean algebra of $\mathbb{B}_1$ generated by $g(P^{N_1}_1\cup A)$. Denote by $\mathbb{\check B}^*$ the
sub-Boolean algebra %$\mathbb{\check B}^*$
of $P^{N_2}_1$ generated   by $P^{N_1}_{1} \cup f^{-1}(A)$.
%Identify  $\bbB^*$ and $\mathbb{\check B}^*$ by extending
% the map $f$ from
%$P^{N_2}_{4,1}- P^{N_1}_{4,1}$ to $A$ defined in the first paragraph  of the proof.

% We also
%denote by $f$  the embedding of $\mathbb{\check B}^*$ into
%$\mathbb{B}_1$ and onto $\mathbb{B}^*$ obtained as follows.
%which is isomorphic to the subalgebra of
%$P^{N_2}_1$ generated by the same set.
Compose $g$ with
 the union of
 the identity on $P^{N_1}_{1}$ with  the map $f$ given
 in the first paragraph  of the proof
 % on
%$P^{N_2}_{4,1}- P^{N_1}_{4,1}$,
using  the
operations of $N_2$ to give a map from  $P^{N_1}_1 \cup (P^{N_2}_{4,1}- P^{N_1}_{4,1})$  into $\mathbb{B}_1$ that takes $\mathbb{\check B}^*$ to $\mathbb{B}^*$. We also denote this map by $f$.

To ease notation, we will suppress $g$ and pretend that $P^{N_1}_1\cup A$ is actually\footnote{Clearly, this could be achieved by choosing a new copy of $\mathbb{B}_1$.}    contained in $\mathbb{B}_1$.
%With this extension we can iden

%We will also
%denote this map by $f$.

{\bf Step 3:}
{\em Construct a Boolean algebra $\mathbb{B}_2$ that is an amalgam of
$P^{M_1}_1$ and $P^{N_2}_1$ over $f(\mathbb{\check B}^*) = \mathbb{ B}^*$ such that the
atoms of $\mathbb{B}_2$ are $P^{\mathbb{B}_1}_{4,1} \cup A$.}
$\mathbb{B}_2$ is a quotient of the pushout $\mathbb{B}'_2$ of $\mathbb{B}_1$
and $P^{N_2}_1$ over the sub-Boolean algebra $\mathbb{B}^*$ of
$\mathbb{B}_1$ generated by $P^{N_1}_1$ and $A$.
%Moreover,
%atoms of $\mathbb{B}_2$ are $P^{\mathbb{B}_1}_{4,1} \cup A$.
The crux
of the proof is the specification of the atoms of  $\mathbb{B}_2$; it allows us to
extend the amalgam of Boolean algebras to an amalgam in $\bK_1$.
% to obtain a Boolean algebra
%$\mathbb{B}_2$ which extends $M_1$ and $f(N_2)$ but has the same
%atoms.
%We consider $P^{N_2}_1$ and
%             $\mathbb{B}_1$ to intersect in $\mathbb{B}^*$.

%{\bf Step 4:} Analyzing the amalgam of Boolean algebras.
 By standard properties of the coproduct (Lemma~\ref{pushout}),
 $\mathbb{B}_1$ and $P^{N_2}_1$ are disjointly embedded over $\mathbb{B}^*$ into their coproduct
 $\mathbb{B}'_2$. We will regard the embedding of $\mathbb{B}_1$ as the identity and
 denote by $f$ the embedding of $P^{N_2}_1$ extending our earlier $f$ mapping
  the sub-Boolean algebra $\mathbb{\check B}^*$ of $P^{N_2}_1$
 %generated by the atoms
 into  $\mathbb{B}_1$. Crucially, while $\mathbb{B}_1$ and $f(P^{N_2}_1)$
are sub-Boolean algebras of  $\mathbb{B}'_2$; they are {\em not}
ideals.

 The atoms of the amalgamation base $\mathbb{B}^*$ remain atoms in
$\mathbb{B}'_2$ as: if $a$ is an atom of $\mathbb{B}^*$ then every
$b_1 \in \mathbb{B}_1$ satisfies $b_1 \wedge a = 0$ or $b_1 \wedge a
= a$ and similarly for $b_2\in P^{N_2}_1$ and therefore also for $b_1
\wedge b_2$; using disjunctive normal form, no element of
 $\mathbb{B}'_2$ contradicts the atomicity of an atom
of $\mathbb{B}^*$.  Recall $N_1 \subseteq_{fr} N_2$ is witnessed by
$(I_2,H_2)$. To guarantee the atoms of $\mathbb{B}_1 \setminus
\mathbb{B}^*$ (i.e. $P^{M_1}_{4,1}- P^{N_1}_{4,1}$) are atoms of
$\mathbb{B}_2$, we divide $\mathbb{B}'_2$ by the
ideal\footnote{Abusing notation, since $\mathbb{B}_1$ is not a
$\tau$-structure, we write $P^{\mathbb{B}_1}_{4,1}$ for the set of
atoms of $\mathbb{B}_1$ and $P^{\mathbb{B}_1}_{4}$ for their finite
joins.}, $\Iscr$, generated by
% .% \sidebar{as
%modified May 3}

 $$ \Iscr_0 =\{a \wedge f(b):a \in P^{\mathbb{B}_1}_{4,1} \setminus P^{N_1}_{4,1}, b\in I_2, a \wedge f(b) < a\}.$$

 \medskip
(*) Since each element of $\Iscr$ is {\em strictly} below
 a finite
join of atoms in $\mathbb{B}'_2$ (actually in ${\mathbb{B}_1}$),
$\Iscr$ is a proper ideal of ${\mathbb{B}_1}$ bounded by elements of
$P^{\mathbb{B}_1}_{4}$; but $\Iscr \cap P^{\mathbb{B}_1}_{4} =
\emptyset$.  Indeed, by freeness of the coproduct, $\Iscr \cap {\mathbb{B}_1} = \emptyset$.
 Note that the subalgebra of $\mathbb{B}'_2$ generated by $f(I_2)$ is a subset of $\mathbb{B}_1$
 so it is disjoint from $\Iscr$.

\medskip
%Thus,  $I_2$ is independent from $P^{M_1}_1$ over $P^{\mathbb{B}_1}_{4}$.
  Let
$\pi$ map $\mathbb{B}'_2$ onto $\mathbb{B}_2 =_{def} \mathbb{B}'_2/\Iscr$.
%Further since each
%generator of $\Iscr$ is strictly below an atom in ${\mathbb{B}_1}$,
%and $I_2$ is independent from $\mathbb{B}_1$ over $\mathbb{B}^*$ any
%element of $\Iscr$ is strictly below an atom in ${\mathbb{B}_1}$ (otherwise a non-trivial identity would hold violating the independence;
By (*), no element of $\mathbb{B}_1 \cup f(I_2)$ is collapsed by the
map $\pi: \mathbb{B}'_2\rightarrow \mathbb{B}_2$.  Thus, $\pi$ is
$1$-$1$ on $\mathbb{B}_1\cup f(P^{N_2}_1)$ and $\mathbb{B}_2$ is a
disjoint amalgamation of the Boolean algebras $\pi(\mathbb{B}_1)$ and
$\pi(f(P^{N_2}_1))$.
% Further
%the atoms of $\mathbb{B}_2$ are those of $\mathbb{B}_1$: No element
%of the ideal is beneath an atom $a \in \mathbb{B}^*$ so those remain
%atoms.  So, w
Since $\mathbb{B}'_2$ is generated by $\mathbb{B}_1\cup
f(P^{N_2}_1)$, without loss of generality, we can assume the preimage
of a potential atom of $\mathbb{B}_2$ has the form $a \wedge f(b)$
where $a \in \mathbb{B}_1 -\mathbb{B}^*$ is an atom of
$\mathbb{B}_1$ and $b\in P^{N_2}_1 -\mathbb{B}^*$. By the
freeness property of coproducts\footnote{$\mathbb{B}'_2$ is freely
generated as a Boolean algebra by (isomorphic copies of) $\mathbb{B}_1$ and $P^{N_2}_1$ over
$\mathbb{B}^*$.}, $\mathbb{B'}_2\models a \wedge f(b) < a$, so $\pi(a
\wedge f(b))= 0$ and $\pi(a)  = a$ is an atom.

% But then the freeness implies
%there are infinitely many elements below $a \wedge b$.
{\bf Step 4:} {\em The actual $\tau$-amalgam.} Now to define the extension
$M_2$, let $P^{M_2}_1 = \mathbb{B}_2$, $P^{M_2}_{4,1} = P^{M_1}_{4,1}
\cup A$; $P^{M_2}_{4}$ is the set of  finite joins of these atoms.
Then, let $P^{M_2}_2 = P^{M_1}_2 \cup P^{N_2}_2$ and the $F^{M_2}_n(c)$
be as in  whichever of $M_l$, $N_2$ in which $c$ lies.  Define $P^{M_2}_0$ to be a set in $1$-$1$
correspondence with $P^{M_2}_{4}$ and call the correspondence
$G^{M_2}_1$. Finally, we must define $R^{M_2}$: for each $b \in
P^{M_2}_1$, let $R({M_2},b) =\{a \in P^{M_2}_0: G^{M_2}_1(a)
\leq^{M_2} b\}$.
%We have shown $M_1 \subseteq_{fr} M_2$.

By Lemma~\ref{potoind}, $I_2$ is independent from $\mathbb{B}_1$ over  $\Iscr$ in $\mathbb{B}'_2$
 and so, by (*),  $\pi(f(I_2))$ is independent from $P^{M_2}_1=\mathbb{B}_2$ over $P^{M_2}_4$ in $M_2$. So $M_1 \subset_{\fr} M_2$ with
 $H_{M_2}(c) = n^*$ for $c\in P^{M_2}_2$.
 $\qed_{\ref{amalfree}}$

\medskip
 Note the $M_0,M_1,M_2,M_3$ in the next argument are
$N_0,M_1,N_2,M_2$ in  Lemma~\ref{amalfree}.

\begin{corollary}\label{corapholds} %[3.2.11 in October, 2016 version]
$(\bK^1_{<\aleph_0}, \subseteq)$ has the disjoint amalgamation property.
\end{corollary}

%\sidebar{apr21: moved end of old proof - must modify}

Proof. We know every member of $\bK^1_{<\aleph_0}$ is $\bK_1$-free over the empty set.  So
the amalgamation becomes a special case of Lemma~\ref{kotofree} when
we add a proof that the amalgam is in $\bK^1_{<\aleph_0}$. We have the following
situation.  $M_0$ is $\bK_1$-free over the minimal model $M_{min}$. That is,
there are $J_0,I_0, H_0$ such that $J_0$ generates $B_{0,n^0_*}$ and
$(J_0 \cup I_0, H_0)$ witness that $M_0$ is a $\bK_1$-free extension of  the
minimal model $M_{min}$. Similarly there are for $i = 1,2$, $J_i,I_i,
H_i$ such that $J_i$ generates $B_{i,n^i_*}$ and $(J_i \cup I_i,
H_i)$  that witness that $M_i$  is a $\bK_1$-free extension of  the minimal
model $M_0$.

%$J_1 \cup I_1$ generates $P^{M_1}_1$ and $J_1 \cup I_1$ is free from
%$P^{M_0}_1$ over $P^{M_1}_4$.
%
%$J_2$ generates $B_{2,n^2_*}$ and is free from $P^{M_0}_1$ over
%$P^{M_2}_4$.
%
%$J_2 \cup I_2$ generates $P^{M_2}_1$ and $J_2 \cup I_2$ is free from
%$P^{M_0}_1$ over $P^{M_2}_4$.
%
Choose  $n_*$ as the maximum of $n^i_*$ for $i<3$; we can assume the
$n^i_*$ for $i<3$ are equal and that $B_{2,n^i_*} \cap B_{0,n^0_*} =
B_{0,n^0_*}$ for $i = 1,2$. Rechoosing $n_*$ by Lemma~\ref{fixbase}
we can assume for all $n \geq n_*$, $ B^1_n
 \cap P^{M_0}_1= B^0_n =B^2_n \cap P^{M_0}_1$.

  Choose $M_3$ by Lemma~\ref{amalfree}.
Let $b^3_{n_*} = b^1_{n_*} \wedge  b^1_{n_*}$. Now let $B_{3,n_*}$ be
the subboolean algebra of $M_3$ generated by $J_0 \cup J_1 \cup J_3$
and for $n \geq n_*$, $B_{3,n}$ be generated by $B_{1,n} \cup
B_{2,n}$.  This is the required witnessing sequence.
$\qed_{\ref{corapholds}}$

\bigskip Since $\bK^1_{<\aleph_0}$ has joint embedding, amalgamation and only
countably many finitely generated models, we construct in the usual
way a generic model. This construction can be rearranged in order
type $\omega$ so by Theorem~\ref{amalfree} and Lemma~\ref{frmodtrans}
the generic is $\bK_1$-free.

\begin{corollary}\label{getgen}  There is a countable generic model $M$ for $\bK_0$. We denote its
Scott sentence by $\phi_M$. Moreover $M$ is $\bK_1$-free.
\end{corollary}

%In fact, since there are only countably many models in $\bK_0$, we
%can expand the vocabulary by adding predicates for appropriate
%$L_{\omega_1,\omega}$-definable relations so that the models of
%$\bK_2$ are the atomic models of the first order theory of the $M$ in
%this vocabulary.

Aligning our notation with earlier sections of the paper we note the
models of $\phi_M$ are rich in the sense defined there.

\begin{definition}\label{richname} We say a model $N$ in $\bK_1$ is rich  if for any $N_1, N_2 \in \bK^1_{<\aleph_0}$ with $N_1 \subseteq N_2$ and $N_1 \subseteq M$, there
is an embedding of $N_2$ into $N$ over $N_1$.  We denote the class of
rich models in $\bK_1$ as $\bK_2$ or $\Krich$.
\end{definition}

Lemma~\ref{amalfree} finds a $\bK_1$-free extension of each $\bK_1$-free model in $\bK_1$; more strongly:

\begin{corollary}\label{getfreeext} Let $M_1$ be $\bK_1$-free.
% \begin{enumerate}
% \item
 There exists an $M_2 \in \bK_2$ which is a proper $\bK_1$-free extension of $M_1$.
 % \item We can choose $M_2\in \bK_2$.
%\end{enumerate}
\end{corollary}

Proof. %For 1) embed the minimal model $M_0$ (Lemma~\ref{minmod}) into
%$M_1$; let $N_1$ be an extension of $M_0$ by adding one more $d$ to
% $P^{M_0}_2$ to form $P^{N_1}_2$ and setting $R({N_1},F_n(d))
%=\emptyset$. Amalgamate $M_1$ and $N_1$ over $M_0$ by
%Lemma~\ref{amalfree}. Note the result $M_2$ is a free extension of
%$M_1$.
Iterate Corollary~\ref{amalfree} as in Corollary~\ref{buildrich} to
obtain a rich model; note that $\bK_1$-freeness is preserved at each stage.
$\qed_{\ref{getfreeext}}$

\medskip

The crucial distinction from Corollary~\ref{buildrich} is that here
we extend only `$\bK_1$-free models' in $\bK_1$ to $\bK_2$.
 While {\em this} construction applied to models in  $\bK_2$ will necessarily increase $P_0$ (case 2 of Construction~\ref{detzfc+}), we can find extensions
in $\bK_1$ which do not extend $P_0$ or $P_1$ but only $P_2$  (case 4 of Construction~\ref{detzfc+}).

For the construction in Section~\ref{zfc+case} we require two
crucial properties of the generic model.

%Now we axiomatize $\bK_2$ and point out some distinctions between the
%theory of $\bK_2$ from that of $\bK_1$.  Note that, in particular,
%the ancillary $B_n$ used to characterize models in $\bK_1$ are not
%used in discussing $\bK_2$. It is easy to check using
%Lemma~\ref{nonoise} the following properties are true in the generic.

%\begin{lemma}\label{K2ax} The Scott sentence $\psi$ of the generic model (and thus
%$\bK_2$ is determined by the following properties that are
%expressible in $L_{\omega_1,\omega}$.

\begin{lemma}\label{nonoise}
If $M$ is the generic model then
\begin{enumerate}[i]\item
if $b_1 \neq
b_2$ are in $P_1^M - P_4^M$ then $R(M,b_1) \neq R(M,b_2)$, i.e.  the map
$f$ from Definition~\ref{defK1}.3 is injective.
 \item For any $a \in P^M_0$, $b \in P^M_1$, $M \models R(a,b) \vee
   R(a,b^-)$. Indeed, $P_1^M$ is an atomic Boolean algebra.
%\begin{enumerate}
%\item $M\models T$, the first order consequence of $\psi$.
%\item $P^M_1/P^M_4$ is an atomless Boolean algebra.
%\item $P^M_{4,1}$ is infinite.
%\item $P^M_4$ is a distributive lattice\footnote{With $b^*_n$ as
%    the `1', $P^M_4$ is a Boolean algebra but it is {\em not} a
%    Boolean {\bf sub}algebra of $P^M_1$.} consisting of the
%    collection of finite joins of  elements of $P^M_{4,1}$ plus
%    the least upper bound $b^M_*$ of $P^M_{4,1}$ (This element is
%    first order definable using the predicates $P^M_{4}$ and
%    $P^M_{4,1}$.)
%\item $P^M_1$ is an atomic Boolean algebra: every non-zero
%    element is above an atom.
%    \item The sets $R(M,F^M_n(c))$ with $F^M_n(c) \not \in P^M_4$
%        and $c  \in P_2(M)$ are independent in the algebra of
%        sets on $P^M_0$.
\item For each $b \in P^M_1 - P^M_4$,
      $R^M(M,b)$ is infinite and coinfinite.
\end{enumerate}
%%\end{enumerate}

\end{lemma}

Proof. For i) fix a finitely generated model $M_0$ containing
$b_1,b_2$; there is a finitely generated extension $M_1$ in $\bK^1_{<\aleph_0}$
by adding $a \in P_0^{M_1}$ with $R^{M_1}(a,b_1) \wedge \neg
R^{M_1}(a,b_2)$. This shows the injectivity;
the other conditions are similar. $\qed_{\ref{nonoise}}$

%the next three
%statements hold by similar arguments. For v) note that $M$ is
%constructed by a chain of $\bK_1$-free extensions and apply
%Lemma~\ref{frmodtrans}.

\medskip
%Lemma~\ref{nonoise}.2 presents an interesting contrast between
%$\bK_0$ and $\bK_2$.  If $N \in \bK_0$ and $n \geq n_*$, $F^N_n(c)
%\wedge b =0$ for any $b \in P^N_4$ so $R(F^N_n(c),N) = \emptyset$.
%But as $n_*$ increases, this condition is eventually violated for
%every $n$ and the contrary  (2) of Lemma~\ref{nonoise}) holds in the limit.

%\sidebar{right theorem??? we are arguing from richness  -the crux is
%qe.  We need these properties in the final proof.}

%\sidebar{Do we have to worry about separability?}
\begin{lemma}\label{compsent}  If $M,N \in \bK_2$,
$M\equiv_{\infty,\omega} N$ so they satisfy the Scott sentence
$\Phi_M$. Moreover, if $M \subset N$ and are both in $\bK_2$,
$M\prec_{\infty,\omega} N$.
\end{lemma}
 Proof: Suppose $M$
and $N$ are in $\bK_2$. We define a back-and-forth between $M$ and
$N$ for $\abar \in M^n$, $\bbar \in N^n$ by $\abar \equiv \bbar$ if
they realize the same first order type over the $\emptyset$ with
respect to $T$. % i.e. the subalgebras they generate are isomorphic.
Fix such $\abar \equiv \bbar$ and choose $c \in M$. The interest is
when $c$ is not in $A =\acl(\abar)$,
% the $\tau$-closure\footnote{The functions
%$H_1$ and $g_n$ in Definition~\ref{defK1} guarantee, e.g. that if $c
%\in P^M_{4,n} - A$ for some $n$, that $c \in P^{A'}_{4,n}$ by adding
%new atoms.} of $\abar$.
If $c \in P^M_1 - A$, let $A_1 = \langle A
\cup \{c\} \rangle_M$. Since $M \in \bK_1$, $A_1 \in \bK_1$.  Now let
$B =\langle \bbar \rangle_N$ that is equivalent to $A$. By richness
there exists $B_1$ isomorphic to $A_1$ with $B \subset B_1 \subset N$.

If $M \subset N$ and both are in $\bK_2$, then $\acl_M(\abar) = \acl_N(\abar)$ for $\abar \in M$; this yields the moreover.
$\qed_{\ref{compsent}}$
%Since the back and forth is determined by quantifier-free types
%%\footnote{More precisely, this becomes true if add,
% we
%have the `moreover'.
%$a'_1, \ldots a'_{k-1} \in P^M_2$ such that $A'$
%the $\tau$-closure of $\abar'$ contains   $A \cup \{c\}$  and the
%$\{F_n(a'_i): n<\omega, i<k\}$ are independent from the $\{F_n(a_i):
%c<\omega, i<n\}$. ????  But then it is easy to choose a similar
%$\bbar' \in N$; the cases $c \in P^M_0, P^M_4, P^M_2$ are analogous
%but easier.

\medskip

This completes our description of the class $\bK_2$ of rich models
and its Scott sentence.
At this point we  show any $\bK_1$-free-member of $\bK_2$ has a proper $\bK_1$-free-extension in $\bK_1$.  In case 2 of Construction~\ref{detzfc+}, we apply Corollary~\ref{getfreeext} to regain a member of $\bK_2$.

%\sidebar{The notion of free extension should either be eliminated or
%modified. Skip to Definition~\ref{nearfree}.}

%  \sidebar{Saharon states this  (9A aprril 6
%scan) for finitely generated $N \in \bK_1$;
% but we need it for arbitrary.}
%
%\begin{lemma} If $M \subset_{fr} N$ and both are in $\bK_1$ then
%\begin{enumerate}
%\item there is an $N_1$ which is a free extension of $N$,
%    $P^{N}_2 =  P^{N_1}_2$, $P^{N_1}_{4,1} = P^{N}_{4,1} \cup
%    \{a\} $.
%
%{\bf above problematic}
%
%
%    \end{enumerate}
%\end{lemma}
%    Proof. 1) Let $P^{N_1}_{4}$ be the collection of element
%    $\{G_1(a)\} \wedge b$ for $b \in P^{N}_{4,1}$
%

%\sidebar{jb adapted from page
%13A of scan A, Apr. 5, 2017.  I don't have any difficulty with this argument but on page 21, the application goes wrong.}

\begin{lemma}\label{nameu} If $M \in \bK_2$, there is an $N$ such that  $M \subset_{fr} N$, both are in $\bK_1$,
  $P^{N}_2 =  P^{M}_2$, $P^{N}_0 =  P^{M}_0$, and
    $P^{N}_{1}$ is generated by $P^{M}_{1} \cup \{b\} $ and $b \in N'$ with  $N \prec N'$.  Moreover given $u
    \subseteq P^{M}_0$, we can require $R(N,b) = u$ and $b$ is free
    from $P^{M}_{1}$ over $P^{N}_{4}$.
     Finally, if $M$ is $\bK_1$-free then so is $N$.
    \end{lemma}

    Proof:
    %For each finite $A$ contained in $P^M_0 \cup P^M_1$, we define a
%    structure $N_A$.
    Let $p(x)$ be the type of an element satisfying $P_1(x) \wedge \neg P_4(x)$:
    $$\{x \geq G_1(a)\colon a \in u\} \cup \{G_1(a) \wedge x =
    0\colon a \in P^N_4 \setminus u\} \cup \{b \wedge \sigma( x)\neq a\colon b
    \in P^{M}_1 \setminus P^M_{4}, a \in P_4^M\},$$
where  $\sigma(x)$ ranges over nontrivial Boolean polynomials.
Each finite subset $q$ of $p$ is satisfied in $M$ because $M\in \bK_2$.
 Thus
there is an elementary extension $N'$ of $M$ where $p$ is realized by
some $b$.  Let $\bbB$ be the boolean subalgebra of $P^{N'}_1$
generated by $P^M_1 \cup \{b\}$. Since $N'$ satisfies the first order properties of $\bK_2$, the atoms of
$M$ are atoms of $\bbB$.

%Oct 23, 2018 -- next paragraph seems irrelevant.
%Moreover, if $d \in \bbB - P_4^{N'}$ then there is an ultrafilter $\Dscr$ on $\bbB$ such
%that for each $c\in P_2^M$, at most finitely many $F_n^M(c) \in \Dscr$. To see
%this, write $d$ as $(d_1 \wedge b ) \vee (d_2 \wedge b)$ for some $d_1,d_2 \in P_1^M$.
%% and by symmetry it suffices to consider $d$ as $(d_1 \wedge b )$.
%Now, since the $F^M_n(c)$ are independent (as sets) in  $P_1^M$, there is an ultrafilter $\Dscr_1$ on $\bbB$ with $(d_1 \wedge b )\in \Dscr_1$ and
% for each $c\in P_2^M$, at most finitely many $F_n^M(c) \in \Dscr_1$. Now extend
% $\Dscr_1$ to an ultrafilter  $\Dscr$ on $\bbB$  with $(d_1 \wedge b )$ and hence $b\in \Dscr$.
% Thus for any $c\in P_2^M$ no element of $\bbB$ realizes $\bigwedge_{n<\omega} F^M_n(c)$. Let $Y$ be the set of atoms of the  atomic Boolean algebra $P_1^{N'}-P_1^M$ that are below an element of $\bbB$.
% And let $\bbB'$ be a subalgebra of  $P_1^{N'}$ generated by $\bbB \cup Y$.
%Now for any $c\in P_2^M$, no element of $\bbB'$ realizes $\bigwedge_{n<\omega} F^M_n(c)$.

Define a $\tau$-structure $N$ with
$P_1^N = \bbB$. Interpret $P_2$ and the
$F_n$ in $N$ as in $M$.  Extend $G_1^M$ and $P_0^M$ so that  $P_0^N = (G_1^N)^{-1}(Y)$.
 The
structure $N$ is well-defined; we must prove it is in $\bK_1$.

Let $\langle (M_i,Z_i)\colon i< |M|, Z_i \subset_\omega Y\rangle$ list the
pairs of finitely generated $M_i \subset M$ in $\bK^1_{<\aleph_0}$ and finite subsets $Z_i$ of $Y$. (The $M_i$ will be repeated.) Let $N_i \subset N$ with
$P^{N_i}_0 = P^{M_1}_0\cup Z_i$, $P^{N_i}_2 = P^{M_1}_2$, and $P^{N_i}_1$ be
the universe of the Boolean subalgebra of $N$ $\tau$-generated by $P^{N_1}_2
\cup \{b\} \cup Z_i$. It is easy to check each $N_i\in \bK^1_{<\aleph_0}$.
Now  $N$ is the direct limit of the finitely generated $\{N_i\colon i< |M|\}$ so it is in  $\bK_1$.

%Then $N_i \in \bK_0$ as $N_i$ is also generated by
%$P^{M_1}_1 \cup \{b_1\}$, where $N \models b_1 =  b-\bigvee Z_i$.

% b-\bigvee\{G_1(a): a
%\in P^{N_1}_0 \cap u\}$.
 Finally $b$ is free from
    from $P^{M}_{1}$ over $P^{N}_{4}$ since no nontrivial unary polynomial $\sigma$ satisfies maps $\sigma(b) \wedge a \in P^N_4$ with $a \in P^M_1-P_4^M$.  The moreover follows by Definition~\ref{deffreesub} from the independence of  $b$.
    $\qed_{\ref{nameu}}$

\subsection{Constructing maximal models in an extension of
ZFC}\label{zfc+case} \setcounter{theorem}{0} We  show that for arbitrarily large
cardinals below a measurable cardinal, assuming a  mild set theoretic
hypothesis described below, $\bK_2$ has  maximal models.
 We begin by defining a pair of set theoretic
notions and some specific notions of maximal model.

%There are two construction depending on the set theoretic hypotheses.
%\sidebar{ Should these give common properties that can be
%summarised?}
%
%
%\subsection{The construction}

\begin{definition}[$\diamond_S$]  Given a cardinal $\kappa$ and a stationary set $ S\subseteq\kappa $,
 $\diamond_S$ is the statement that there is a sequence $\langle A_\alpha: \alpha \in S \rangle $ such that
\begin{enumerate}
\item each $ A_\alpha \subseteq \alpha $; \item for every $ A
    \subseteq \kappa, \{\alpha \in S: A \cap \alpha = A_\alpha\}
    $ is stationary in $\kappa$.
\end{enumerate}
\end{definition}
\begin{definition}[$S$ reflects]
Let $\kappa$ be a regular uncountable cardinal and let S be a
stationary subset of $\kappa$. For $\alpha< \kappa$ with uncountable
cofinality, $S$ reflects at $\alpha$ if $S\cap \alpha$ is stationary
in $\alpha$ . $S$ reflects if it reflects at some $\alpha< \kappa$.
\end{definition}
%\sidebar{Shouldn't it be $P_0$ maximal in next line?}

%Note that if the elements of $S$

%\sidebar{aug 9:Presumably $\bK^{BA}_1$ in the next definition should be
%$\bK_2 = \Krich$ and there are similar mistake spread around.
%
%aug 10: made this change and changed $\subseteq_*$ to $\subseteq$.}

\begin{definition}\begin{enumerate}
\item A model $M \in \bK_2 = \Krich$ is {\em $P_0$-maximal} (for
    $\bK_1$)
    if $M \subseteq N$ and $N \in \bK_2$ ($\in \bK_1$) implies $P_0^M =
    P_0^N$.
\item A model $M \in \bK_2$ is {\em maximal} for $\bK_2$
    if $M \subseteq N$ and $N \in \bK_2$ implies $M =
    N$.
\end{enumerate}
\end{definition}

 Let $S^\lambda_{\aleph_0}$ denote the stationary set
$\{\delta< \lambda: \cf(\delta) = \aleph_0, \delta {\rm\  is \
divisible \ by \ |\delta|}\}$.

%\sidebar{jb jun27: I don't see where the divisibility condition is
%used.}

%\sidebar{Presumably $\widehat{\bK}$ in the next definition should be
%$\bK_2 = \Krich$ and there are similar mistake spread around.}

\bigskip

We now define a crucial notion.

%\sidebar{must connect with the $\gammabar$ -- maybe later}

\begin{definition}[$A$-good defined]\label{gooddef}  Suppose that $N_n \subset_{\fr} N_{n+1}$ for $n< \omega$,  is sequence of models, $\overline{N}$,  in $
\bK_1$.  We say  a sequence $\bbar= \langle b_n\colon n<
    \omega\rangle$ is

\begin{enumerate}
\item
    {\em good for  $\overline{N}$} if
\begin{enumerate}
\item  $P^{N_{n+1}}_2 - P^{N_n}_2$ is infinite;
%\footnote{It would suffice
%if there are no finite $Y\subset P^{N_{n+1}}_1$ and
% $X \subset P^{N_{n+1}}_2$ so that $P^{N_{n+1}}_1$ is
% generated as a Boolean algebra by $P^{N_{n}}_1 \cup Y \cup
% \{F^{N_{n+1}}_m(c): c \in X, m \in \omega\}$};
 \item for each $n$, $b_n \in P^{N_{n+1}}_1$ and
     $\{b_n\}$ is free from $P^{N_{n}}_1$ over
     $P^{N_{n+1}}_4$;

\item if $a \in  P^{N_{i}}_0$, then for all but finitely many
    $n\geq i$, $a \not \in R({N_{n+1}},b_n)$. \label{empint}
\end{enumerate}

\item for $A \subset \bigcup \overline{N}$, $\bbar$ is {\em $A$-good} if each $b_n \in A$.

\item and {\em labeled} if there is a pair $(N^{\bbar},c^{\bbar}) $ = with $N^{\bbar} \in
    \bK_1$ and  $N^{\bbar} \supseteq N_\omega = \bigcup N_n$ such
    that for each $n$, $F^{N^b}_n(c^{\bbar}) = b_n$.  By the definition of $\bK_1$, $\bigcap_n R(N^{\bbar},F^N_n(c^{\bbar})) = \emptyset$.
\end{enumerate}

\end{definition}

Note that for every $c \in N_m \subsetneq N_\omega$, at most finitely many of any good sequence $\langle b_k\colon k<\omega\rangle$ occur in the sequence $F^{N_m}_n(c)$ for $n<\omega$ (as $F^{N_m}_n(c) \in N_m$ and for $k>m$, $b_k \not \in N_m$).
%\sidebar{Is the really plausible? Nothing is said in the proof about $R(N^{\bbar},b_n)$.}

 Any proper $P_0$-extension of a model $M$ induces a  non-principal ultrafilter $A$ on $P^{M}_1$.  Claim~\ref{cl2} is instrumental via case 5 in constructing, for the particular $M$ under consideration,    an ostensibly  non-principal $\aleph_1$-complete ultrafilter on  $\Pscr(P^{M}_0)$ which contradicts that $\lambda$ is not measurable.  See \ref{verzfc+}.

\begin{claim}\label{cl2}   Suppose that for $n< \omega$, $\overline{N} = \langle N_n \subset_{\fr} N_{n+1}\rangle$ are in $
\bK_1$.  For $A \subseteq N_\omega$, if Condition A) holds then so does condition B).
\begin{itemize}
\item [A)]  There is an $A$-good sequence for $\overline{N}$.
\item [B)]   There is a  labeled $A$-good sequence for $\overline{N}$.

%There is a sequence $\bbar= \langle b_n\colon n<
%    \omega\rangle$ such that
%\begin{enumerate}
%\item  $P^{N_{n+1}}_2 - P^{N_n}_2$ is infinite;
%
%or at least there are no finite $Y\subset P^{N_{n+1}}_1$ and
% $X \subset P^{N_{n+1}}_2$ so that $P^{N_{n+1}}_1$ is
% generated as a Boolean algebra by $P^{N_{n}}_1 \cup Y \cup
% \{F^{N_{n+1}}_m(c): c \in X, m \in \omega\}$;
% \item for each $n$, $b_n \in P^{N_{n+1}}_1$ so that
%     $\{b_n\}$ is free from $P^{N_{n}}_1$ over
%     $P^{N_{n+1}}_4$;
%
%\item if $a \in  P^{N_{i}}_1$, then for all but finitely many
%    $n\geq i$, $a \not \in R({N_{n+1}},b_n)$.
%\end{enumerate}
%
%\item [B)] then %for each such $\bbar$,
%there is a pair $(N^{\bbar},c^{\bbar}) $ = with $N^{\bbar} \in
%    \bK_1$ and  $N^{\bbar} \supseteq N_\omega = \bigcup N_n$ such
%    that:
%
%\begin{enumerate} \item  $N^{\bbar}= \bigcup N_n \cup \{c^{\bbar}\}$, $c^{\bbar}\in P^M_2$,
%$c^{\bbar}$ is not in any $N_n$;
%\item $N_n \subset_{\fr} N$ for each $n$;
%\item $F^N_n(c^{\bbar}) = b_n$.
%\end{enumerate}
\end{itemize}
\end{claim}

Proof.  The following construction is for the fixed $A$-good sequence $\bbar$. Let $ N = N_\omega
\bigcup_{n<\omega}N_n$.
% for ease of reading, the superscript $\bbar$ is
%omitted in the proof.
%The construction will be performed for each
%such $\bbar$.
Note that each $P^{N^{\bbar}}_1 = P^N_1$; the extension
$N^{\bbar}$  only adds an element $c$ to $P^{N_\omega}_2$ and
interprets the $F^{N^{\bbar}}_m(c)$.
The difficulty is that while we know each $N_{n+1}$ is $\bK_1$-free over
$N_n$, witnessed by some $(I_n, H_n)$, we don't know $b_n \in I_n$.
We need to find $I'_n$ which  witnesses both
 $N_n \subset_{\fr} N_{n+1}$ and $b_n \in  I'_n$. After this construction we will choose an $N^{\bbar}$ extending $N$ witnessing goodness.

% Then we choose
%$c^*$ and define $N$  with domain $\bigcup N_n \cup {c^*}$ by keeping
%the $\tau$-structure on $\bigcup N_n$, adding $c^*$ to $P^N_2$ and
% setting $F_n^N(c^*) = b_n$ for each $n$.
 To find $I'_n$,
we first find $(X_n,J_n)$ such that:
\begin{enumerate}
\item $X_n \subseteq P^{N_n}_1$ is finite.
\item $J_n \subset I_n$ is countable.
\item If $c \in P^{N_{n+1}}_2 - P^{N_n}_2$ then for sufficiently
    large $m$, $F^{N_{n+1}}_m(c) \not \in J_n$.
\item $b_n \in BA( X_n\cup J_n)$, the {\em Boolean algebra
    generated } by $X_n\cup  J_n$ in $P^{M_{n+1}}_1$.
\end{enumerate}

{\em First step}: First, we construct such an $(X_n,J_n)$. Note that
$b_n$ is in  a subalgebra generated by a finite subset $X_n$ of
$P^{N_n}_1$ and a finite subset $J'_n$ of $I_n$.

Now, by 1a) of Definition~\ref{gooddef} , fix a sequence $\langle c_i\colon i < \omega\rangle$ of
distinct elements of $P^{N_{n+1}}_2 - P^{N_n}_2$. Note that for $i,j<
\omega$ if $n_i
> H_n(c_i)$ and $n_j
> H_n(c_j)$ then $F^{N_{n+1}}_{n_i}(c_i) \neq F^{N_{n+1}}_{n_i}(c_i)$.  Now we can construct
a $J''_n =\{d_{n,k}:k< \omega\}$ from $I_n -J'_n$ by $d_{n,k} =
F^{N_{n+1}}_m(c_k)$ for some $m > H_n(c_k)$. We now have a countably
infinite $J''_n$ contained in $I_n -J'_n$ such that for each $c \in
P^{N_{n+1}}_2 - P^{N_n}_2$ all but finitely many of the
$F^{N_{n+1}}_m(c)$ are in $I_n -(J'_n \cup J''_n)$.   Set $J_n =J'_n
\cup J''_n$.

%\sidebar{Why is the last sentence important?  Can we get away with
%arb large witnesses?}

% since the
%$F^{M_{n+1}}(c)$ are distinct.

%\sidebar{Still need  page 2A.5}
{\em Second step}:  Now apply Lemma~\ref{basisext}\footnote{This  is the crucial application of Lemma~\ref{basisext} which stengthened our notion of independence by getting a standard consequence of exchange, even though exchange fails here.} to find $J_n^*$
with $J_n^*$ independent from $P^{M_n}_1$ over ${P^{M_{n+1}}_4}$ such
that $\langle J_n^* \cup P^{M_{n+1}}_4 \rangle_{P^{M_{n+1}}_1} =
\langle J_n \cup P^{M_{n+1}}_4 \rangle_{P^{M_{n+1}}_1}$ but $b_n \in
J_n^*$. Now, $I'_n$ can be taken as $(I_n - J_n) \cup J_n^*$. To
ensure that $N_n \subseteq_{fr} N_{n+1}$ with basis $I'_n$, replace
$H_N(c_n)$ by $H_{N_n}(c_n) +r_n$ where (by Definition~\ref{k0}) some
$r_n$ bounds the number of $m$ such that $F^{N_n}_m(c_n) \in \langle
J_n \rangle_{P^{M_{n+1}}_1}$.

%\sidebar{But we can't iterate as we might lose some each time??? Does
%it only depend on $\gammabar$ and is that enough?}

Having found an appropriate basis for $N = \bigcup N_n$, we extend
$N$ to $N^{\bbar}$ by adding an element $c^{\bbar}$ to $P^{N_2}_2$
and defining $F^{N^{\bbar}}_n(c^{\bbar}) = b_n$.
%Easily,
%$N^{\bbar}\in \bK_1$, as $\bigcup N_n$ is $\bK_1$-free and
%the additional $F^N_n(c^*)$
%just picks out independent elements;
The sentence immediately before Claim~\ref{cl2} guarantees that $N^{\bbar}$ is $\bK_1$-free;
set $H^{N^{\bbar}}(c^*) = 0$; thus, $N^{\bbar}\in \bK_1$.
Since the same $b_n$ were used, it is clear the labeled sequence is $A$-good.
(Note that there is no requirement that $m,n < \omega$, $c\in P_2^{M_0}$, $d\in P_2^{M_1}$ imply $F_n^{M_1}(c) \neq F_m^{M_1}(d)$; we only require that there be only finitely many such conflicts.)
$\qed_{\ref{cl2}}$

\medskip

We now state precisely the main theorem.

\begin{theorem}\label{zfc+thm} Fix $\bK^1_{<\aleph_0}, \bK_1 = \hat\bK^1_{<\aleph_0}$, and $\bK_2 = \Krich$ as in Definitions~\ref{k0}, \ref{k1def} and
\ref{richname}. There is a $P_0$-maximal for $\bK_2$ model $M \in
\bK_2$ of card $\lambda$ if there is no measurable cardinal $\rho$
with $\rho \leq \lambda$, $\lambda = \lambda^{< \lambda}$, and there
is an $S \subseteq S^{\lambda}_{\aleph_0}$, that is stationary
non-reflecting, and $\diamond_S$ holds.
\end{theorem}

Under $V=L$, the hypotheses are clearly consistent and imply there
are arbitrarily large maximal models of $\Krich$ in $L$. When a
measurable cardinal exists, the consistency of the conditions can be
established by forcing; see the article by Cummings in the Handbook
of Set Theory \cite{cummings} or by considering the inner model of a
measurable $L[D]$ where is $D$ is a normal ultrafilter on $\mu$.

The argument for Theorem~\ref{zfc+thm} will have three parts. First,
we describe the requirements on a construction of a rich model; then we
carry out the construction. Finally, we show the model  constructed
is $P_0$-maximal when $\lambda$ is below the first measurable and satisfies the other conditions of Theorem~\ref{zfc+thm}.

\begin{construction} [Requirements]\label{goalzfc+}
{\rm
\medskip

Fix $\lambda$ satisfying the cardinal requirements in Theorem~\ref{zfc+thm}.
List $[\lambda]^{<\lambda}$, the subsets of $\lambda$ with less than $\lambda$ elements, as  $\langle U_\alpha:\alpha< \lambda\rangle$
 so that each subset is enumerated $\lambda$
times and $U_\alpha \subseteq \alpha$. Since the set of ordinals
$\alpha < \kappa$ such that $|\alpha|$ divides $\alpha$ is a cub for
any $\kappa$, without loss of generality, each $\alpha \in S$ is a
limit ordinal and is divided
by $|\alpha|$. %\sidebar{says $\alpha$ --What's meant.}
Let $\overline{A}^* = \langle A^*_\delta: \delta \in S\rangle$ be a
$\diamond_S$-sequence.

%Recall that for a given cardinal number $\kappa$ and a stationary set
%$S\subseteq\kappa$ , the statement $\diamond_S$ asserts that there is
%a sequence $\langle A_\alpha: \alpha \in S \rangle$  such that
%\begin{enumerate}\item each  $A_\alpha \subseteq \alpha$ \item for every  $A \subseteq
%\kappa, \{\alpha \in S: A \cap \alpha = A_\alpha\}$  is
%stationary in $\kappa$.
%\end{enumerate}

\medskip
We will choose $M_{\alpha}$ for $\alpha< \lambda$ by induction to
satisfy the following conditions. (Since the universe of $M$ is a
subset of $\lambda$, its elements are ordinals so we may talk about
their order although the order relation is not in $\tau$.)

%\sidebar{yellow pages didn't specify the universe in 1}

%\sidebar{Presumably $\widehat{\bK}$ in the next definition should be
%$\bK_2 = \Krich$ and there are similar mistake spread around.}

\begin{enumerate}
\item $M_0$
    is isomorphic to the minimal model of $\bK_1$. For $\alpha>1$, $M_\alpha \in \bK_2$ %\in \Krich$
has universe an ordinal between $\alpha$ and $\lambda$.
% Each
%$M_\alpha$ is in $\bK_1$.
\item $\langle M_\beta: \beta < \alpha\rangle$ is $\subseteq$-
    continuous. %\sidebar{Should this be $\subseteq_*$?}

\item If $\beta \in  \alpha - S$ then $M_\alpha$ is $\bK_1$-free over
    $M_\beta$,  {\em and} $M_\alpha \in \bK_2 =\Krich$.

    %\sidebar{Where is condition iii) used?}
    \item If $\alpha = \beta+2$ and $U_\beta \subseteq
        P^{M_\beta}_0$ then there is  a $b_\beta\in
        P^{M_\alpha}_1$ such that $R(M_\alpha, b_\beta) \cap
        M_{\beta +1} = U_\beta$ and in the Boolean algebra
        $P^{M_\alpha}_1$, $\{b_\beta\}$ is free from
        $P^{M_{\beta+1}}_1$ modulo $P_4^{M_\alpha}$.  Moreover
        $P^{M_\alpha}_2 -P^{M_\beta}_2$ is infinite.

%\sidebar{This remark is nonsense $R(N,b)$ extends $R(M,b)$ if $b
%contained in N so independent elements can have same
%intersection. This freeness is clearly impossible if some $b \in
%P^{M_{\beta +1}}$ satisfies $R(M_{\beta +1}, b) = U_\beta$. Do we
%only make this requirement if $U_\alpha$ is not named?}

\item If $\delta \in S$ and $\alpha = \delta+1$ then A) implies
    B), where:

     \begin{itemize}
 \item [A)] there is an $A$-good sequence $\gammaover =
     \langle \gamma_{ \delta,n}, b_{\delta,n}:n< \omega
     \rangle$, where the $\gamma_{ \delta,n}$ are increasing
     with $n$ and not in $S$   such that the $\langle b_{\delta,n}:n< \omega
     \rangle$ are good for the $M_{\gamma_{\delta,n}}$.

     \item[B)]  there is a labeled $A$-good sequence $\hat \gammaover =
     \langle \hat  \gamma_{ \delta,n}, \hat  b_{\delta,n}:n< \omega
     \rangle$, for $\langle M_{\gamma_{\delta,n}}:n< \omega\rangle$ with $c \in M_{\delta+1}$.

 %\begin{itemize}
% \item [A)] there is an increasing sequence $\gammaover =
%     \langle \gamma_{ \delta,n}, b_{\delta,n}:n< \omega
%     \rangle$, where the $\gamma_{ \delta,n}$ are increasing
%     with $n$ and not in $S$   satisfying:
%
%
%
%\begin{enumerate}
% \item $\gamma_{ \delta,n} < \gamma_{\delta,n+1} <
%     \delta$ with $\sup_n \gamma_{ \delta,n} = \delta$;
% \item $ b_{ \delta,n} \in P_1^{M_{\gamma_{\delta,n+1}}}
%     \cap A^*_\delta$ and so $ b_{ \delta,n}\in
%     P_1^{M_{\delta}}$ ;
%
%\item for each $n$, $\{ b_{ \delta,n}$
%% :n<\omega\} $
% is independent from $P^{M_{\gamma_n }}_1 $ over $
%    P^{M_{\delta }}_4$;
% \item if $a \in P_0^{M_\delta}$ then  for all but
%     finitely many $n$,  $\neg R(a,b_{\delta,n})$.
%
%%  old   if $\beta< \delta$ then  for all but finitely man
%%     $n$, $\neg R(\beta,b_{\delta,n})$.
% \end{enumerate}
% \item[B)] For some $\hat {\gammaover}= \langle \hat \gamma_{
%     \delta,n},\hat b^{\delta}_n:n< \omega \rangle$
%     satisfying A (a)-(d), there is a $\hat c_\delta \in
%     P^{M_{\delta+1}}_2$ such that for each $n$,
%     $F^{M_{\delta+1}}_n(\hat c_\delta ) = \hat
%     b_{\delta,n}$.
 \end{itemize}
 \end{enumerate}
}
\end{construction} %\ref{goalzfc+}
\bigskip

%\sidebar{Need to guarantee that the chain is at least cofinally in $\K_2$ and probably on a club.}

%Note that each $M_i$ for $i<\lambda$ is $\bK_1$-free; but, as we will see,
%$M_\lambda$ may not be.
%For $i \not \in S$, this is requirement 3;  if  $i =\delta \in S$, we can guarantee the $M_\delta$ is $\bK_1$-free
%by choosing a sequence
% of successor ordinals $\gamma_n$ with limit $\delta$. By induction each $M_{\gamma_n}$ is $\bK_1$-free over $M_{\gamma_{n-1}}$
% so by Lemma~\ref{frmodtrans}.2 so is $M_\delta$.

\begin{remark} {\rm Condition 5 asserts that for any $A \subseteq \bigcup_{n< \omega} M_{\gamma_{\delta,n}}$: if there is an $A$-good sequence
then there is a labeled $A$-good sequence.  In the proof of Claim~\ref{cl2} we, in fact, took the same sequence
so the `$A$' is preserved automatically. But for each $\delta$ we construct only one pair of a $c$ labeling
a sequence $b_{\delta,n}$.  We fix the relevant $A$ for application in the first  paragraph
of \ref{verzfc+}; it will be an ultrafilter on $P^M_1$ induced by a proper extension.}
\end{remark}

 We now carry out the inductive construction.

\begin{construction}{Details}\label{detzfc+}

%By Lemma~\ref{getfreeext}.2 we can insist each model constructed
%below is in $\Krich$.

{\rm
\medskip
%\sidebar{Maybe $M_0$ is the prime model???}
\noindent {\em Case 1:} $\alpha=0$. Let $M_0$ be the minimal model
from Lemma~\ref{minmod}.  The generic can be taken as $M_1$.

\noindent {\em Case 2:} $\alpha=\beta +1$ and $\beta \not \in S$. If
$\beta$ is a limit we only have to choose, by Lemma~\ref{getfreeext},
$M_\alpha$ to be a $\bK_1$-free extension of $M_\beta$ in $\bK_2$. If
$\beta$ is a successor, there is an additional difficulty.  If $U_\beta
\subset P_0^{M_\beta}$; we must choose $b_\beta$ to satisfy condition
4) and with $M_{\alpha} \in \bK_2$. For this, apply
 Lemma~\ref{nameu} with $U_\beta$ as $U$ and $M_{\beta+1}$ as $M$ to construct $N$ and $b_\beta$. Now iterate Corollary~\ref{getfreeext} $|M_{\beta+1}|$ times to obtain $M_\alpha\in \bK_2$.  This iteration also ensures $P^{M_\alpha}_2 -P^{M_\beta}_2$ is infinite.

% \sidebar{Supposedly the details in the last two sentences solve:
%
% A1Problem:  Applying Corollary~\ref{getfreeext}, $R(M_\alpha, b_\beta)$ is no longer equal to $U_\beta$ but bigger. It must be to guarantee independence but then we haven't met the conditions of the construction.}

\noindent {\em Case 3:} $\alpha=\delta$, a limit ordinal that is not
in $S$. Set $M_\delta = \bigcup_{\gamma<\delta} M_\gamma$. We must
prove that if $\beta \in \delta\setminus S$ then $M_\delta$ is $\bK_1$-free
over $M_\beta$. Since $S$ does not reflect there exists an increasing
continuous sequence $\langle \alpha_i:i< \cf(\delta)\rangle$ of
ordinals less than $\delta$, which are not in $S$ and with $\alpha_0
= \beta$. By the induction hypothesis, since $\alpha_j \not \in S$, for each $i<j< \cf(\delta)$,
$M_{\alpha_j}$ is $\bK_1$-free over $M_{\alpha_i}$. And  by Lemma~\ref{frmodtrans},
$M_\delta$ is $\bK_1$-free over $M_\beta$ as required.

\noindent {\em Case 4a:} $\alpha=\delta+1$, $\delta \in S$, and
clause 5A fails. This is just as in case 2.

\noindent {\em Case 4b:} $\alpha=\delta+1$, $\delta \in S$, but
clause 5A holds.

 So, suppose $\<
M_\beta,b_\beta\>$ for $\beta< \delta$
%\sidebar{ (only for second
%successors)}
have been defined. If there exists $\gammaover$ as in condition 5A)
of Construction~\ref{goalzfc+} we must construct $\hat {\gammaover}=
\langle \hat \gamma_{
     \delta,n},\hat b^{\delta}_n:n< \omega \rangle$ and $\hat c^\delta$ to satisfy condition 5B).
     Take any   $\langle \gamma_{ \delta,n}, b_{\delta,n}:n< \omega
     \rangle$ satisfying 5A.  Let the $M_{\gamma_n}$ be the $N_n$ from  Claim~\ref{cl2} and by that claim,
  %  Hypothesis 1) of Claim~\ref{cl2} is guaranteed by condition 4) (in case 3 of the construction).  Hypotheses 2)
%and 3) of Claim~\ref{cl2}  are conditions c) and d) of 5A) from
%Construction~\ref{goalzfc+}.
choose $M_{\delta+1}$, $\hat c_\delta
\in
     P^{M_{\delta+1}}_2$ such that for each $n$,
     $F^{M_{\delta+1}}_n(\hat c_\delta ) = \hat
     b_{\delta,n}$.
   %   We
%deduce the result from Claim~\ref{cl2}.

\noindent {\em Case 5:}
 Recall that $\delta$
is divisible by $|\delta|$ so we can choose the $\gamma_n$ so that
$\gamma_{n+1} \geq \gamma_n + \omega$ and each $\gamma_n$ is not in $S$.  So, by iterating as in
Corollary~\ref{getfreeext}, $P^{M_{\gamma_{n+1}}}_2 -
P^{M_{\gamma_{n+1}}}_2$ is infinite.
Moreover, again since each $\gamma_n$ is not in $S$,  $M_{\gamma_{n+1}}$ is $\bK_1$-free over $M_{\gamma_{n}}$ so by Lemma~\ref{frmodtrans}, $M_\delta$ is $\bK_1$-free.

%\sidebar{Must add a condition guaranteeing each $P^{M_{\alpha+1}}_2 -
%P^{M_{\alpha}}_2$ is almost always infinite -- not at successors of
%$\delta$ in $S$.}

}

%\sidebar{Add transition.}

%\sidebar{The next claim was labeled 2.5 in 10/16/15.}

\medskip

This completes the construction. We fix the domain of $M$ as the $\lambda$ chosen for Construction~\ref{goalzfc+}.
\end{construction}

\bigskip

\begin{claim}\label{inrich} The structure  $M = \bigcup_{i<\lambda}M_i \in \bK_2$.
\end{claim}

Proof. Since we required the extension to be in $\bK_2 =\Krich$ in
requirement 3 of Construction~\ref{goalzfc+}, for cofinally many $i$,
$M_i \in \bK_2$. By Lemma~\ref{compsent}, they are
$\infty,\omega$-elementary extensions. Hence $M \in \bK_2$.
$\qed_{\ref{inrich}}$

\begin{construction}{Verification that the construction suffices}\label{verzfc+}
{\rm
\medskip

 Now we now show that $M$ is
$P_0$-maximal for  $\bK_2$. Suppose for contradiction there
exists $N$ in $\bK_2$  extending $M$ %($\bK_2$)
such that $P^N_0 \supsetneq P^M_0$. Choose $a^*\in P^N_0 - P^M_0$.
Let $$A = \{b \in  P^M_1\colon R^N(a^*,b)\}.$$ Then, by
Lemma~\ref{nonoise}.ii, for every $a \in P^N_0 $, in particular
$a^*$ and every $b \in P^N_1$ (and so every $b \in P^M_1$) either
$R^N(a^*,b)$ or $R^N(a^*,b^-)$. Thus, the subset $A$ of $P^M_1$ is a
non-principal ultrafilter of the Boolean algebra $P^M_1$. For, if $A$ is
principal, it is generated by some atom $b_0 \in P^{M_1}_4$. Then
$b_0$ must be in $P_{4,1}$ and so $\neg R^N(a^*,b_0)$, contrary to
the hypothesis that $a \in P_0$.  We will show  that $A$ induces an $\aleph_1$-complete
ultrafilter on $\Pscr(P^{M_{\alpha^*}}_0)$ for some
$\alpha^*<\lambda$. But this contradicts that $\lambda$ is below the first
measurable.

  Recall that the $A_\delta$
are the diamond sequence fixed in requirement \ref{goalzfc+} and that $S \subseteq S^\lambda_{\aleph_0}$. Note
$$S_A = \{\delta \in S: M_\delta { \rm\  has\ universe\ } \delta \
\&\ A_\delta = A \cap \delta\}$$ is a stationary subset of
$\lambda$. In the construction, we chose $b_\alpha$ for $\alpha <
\lambda$ which satisfied requirement 4 of
Construction~\ref{goalzfc+}. Note

%$\begin{notation}\label{defC}
$$C = \{\delta< \lambda: \delta
\hbox{\rm\ limit\ } \ \&\ \alpha< \delta \rightarrow b_\alpha<
\delta\}$$
 is a club on $\lambda$. %$\end{notation}

%\begin{remark} If we required that the sequence $M_n$ was elementary
%in $L_{\infty,\omega}$ or even in a countable fragment containing
%$\phi$, the proof of non-maximality would be much simpler. If a
%sequence $R(M_\delta,b_n)$ had empty intersection, in some $R(M,b_n)$
%would also have empty intersection since this expressed by an
%$L_{\omega_1,\omega}$-sentence $\psi$.  But, we work only with
%submodels and the role of $P_2$ and the $F_n$ is essentially to add a
%Skolem sentence for $\psi$.
%\end{remark}

%
%\sidebar{ The previous remark was motivated by the following line of
%thought.  What is the notion of strong submodel for $\widehat{\bK}$? It seems
%to be just substructure??? But why, particularly if the sentence is
%complete $\prec_{\omega_1,\omega}$ is the natural notion for which we
%want maximality.  The next claim assumes $\prec_{\omega_1,\omega}$ is
%strong substructure.
%
%
%
%
%\begin{claim}\label{cl} The existence of $a$ as above is incompatible
%with \ref{diaconst}5.ii and iv)  ( $*\alpha$  e $\alpha$ 3 and 5 on
%2.12 of yellow pages).
%\end{claim}
%
%Proof.  Item ii) says $b_{ \delta,n} \in P^{M_{\gamma_n +1}} \cap
%     A^*_\delta$, so clearly $R^{N}(a,b_{ \delta,n})$
%
%But item 4 says, that the $b_{ \delta,n}$ have empty intersection in
%$M_\delta$; then they also so have empty intersection in any
%$L_{\omega_1,\omega}$-elementary extension.
%
%$\qed_{\ref{cl} }$
%
%It is unclear to me whether the argument below is aimed at ensuring
%that just submodel is enough or whether it also uses an infinitary
%condition.
%
%
%}

There are two cases. We will show the first is impossible and the
second implies $\lambda$ is measurable, contrary to hypothesis. So the
construction yields a $P_0$-maximal model in $\bK_2$.

{\bf Case i)}: For every $\alpha< \lambda$ there is a $ b \in P^M_1
\cap A$ such that $R(M,b)$ is disjoint from $\alpha$ and $\{b\}$ is
independent from $P^{M_\alpha}_1$ over  $P^M_4$.
%For $\alpha< \lambda$, choose $b_\alpha$ witnessing the condition;
%each $b_{\alpha} \in A$ and the sequence is independent in the strong
%sense of Condition  \ref{diaconst}5.ii.

Choose $\delta^* \in S_A \cap C$. %\cap C_1$
%($C_1$ is defined in
%Fact~\ref{someall}.)
Since $\delta^*$ has cofinality $\omega$ we can choose a sequence
$\langle \hat \gamma^{\delta^*_n}: n< \omega\rangle$ such that each is a successor (so not in
$S$), and, as we are in case i), with
$b_{\hat \gamma^{\delta^*}_n}<\hat \gamma^{\delta^*}_{n+1}$. Since condition 5B)
holds there are $\hat c_\delta \in
     P^{M_{\delta+1}}_2$ such that for each $n$,
     $F^{M_{\delta+1}}_n(\hat c_\delta ) = b_{\hat \gamma^{\delta^*}_n}$. Since $M_{\delta+1} \in \bK_1$, by clause 8 of Definition~\ref{defK1}, $M_{\delta+1}
\models \neg (\exists x)\bigwedge_n R(x, F_n(c_\delta^*))$. This
contradicts that we chose $b_{\hat \gamma^{\delta^*}_{n}} \in A$, since by
the definition of $A$, for each $n< \omega$,
$R^N(a^*,b_{\hat \gamma^{\delta^*}_{n}})$ holds.

\medskip
{\bf case ii)}  For some $\alpha^*$, there is no such $b$. That is,
if $b \in P_1^M$ is independent from $P_1^{M_{\alpha^*}}$ over
$P_4^{M}$ and $R(M,b)$ is disjoint from $\alpha^*$ then $\neg
R(a^*,b)$.
From the list of elements of $[\lambda]^{<\lambda}$ at the beginning of Construction~\ref{goalzfc+}, we consider the subsequence $\langle
v_\gamma\colon \gamma < \lambda\rangle$ enumerating
$\Pscr(P^{M_{\alpha^*}}_0)$; recall each element appears $\lambda$
times in the list.

%\sidebar{double subscript version below}

We now choose inductively by requirement 4 of Construction~\ref{goalzfc+} and
Lemma~\ref{nameu} a subsequence\footnote{   For local intelligibility (and at the risk of
 global confusion) we use indices $b_{\gamma}$ and $M_{\gamma}$ rather than
 $b_{\alpha_\gamma}$ and $M_{\alpha_\gamma}$ that would keep more precise track of the subsequence fact.} $b_\gamma$ %(technically $b_{\alpha_\gamma}$
 of the $b_\alpha \in
P^M_1$
 and $M_\gamma$  such that
 $b_\gamma \in P^{M_{\gamma+1}}_1$ and  $R(M_\gamma,b_\gamma) \cap P^{M_{\alpha^*}}_0 =
 v_\gamma= R(M,b_\gamma)$
%  and so $R(M,b_\gamma) \cap P^{M_{\alpha^*}}_0 =
% v_\gamma$
%\sidebar{jbsep5: Replaced typo in SS copied in my earlier versions
%
% such that $R(M,b_\gamma) \cap P^{M_{\alpha^*}}_0 =
% v_\gamma$ }
%
%
%
 %(moreover $R(M,b_\gamma) \cap P^{M_{\gamma+1}}_0  =
 %v_\gamma$)
 and $\langle b_\beta\colon \beta \leq \gamma\rangle$ is
 independent from  $P^{M_{\alpha^*}}_1$ over $P^M_4$. % with $P^M_1$.
 In
 particular, $b_\beta$ is  independent from $P^{M_{\beta}}_1$ over $P^{M_{\beta+1}}_4$ and
 so by Remark~\ref{monnote} over $P^M_4$.
We claim   that if $\gamma_1 < \gamma_2 \wedge v_{\gamma_1} =
v_{\gamma_2}$ then $R^N(a^*,b_{\gamma_1}) \leftrightarrow
 R^N(a^*,b_{\gamma_2})$. For this, let $b' = b_{\gamma_1} \triangle b_{\gamma_2}$.
Then $R(M,b')  \cap P^{M_{\alpha^*}}_0 = \emptyset$ so by the case
choice, $\neg R(a^*,b')$.  But, as required, $\neg R(a^*,b')$ implies
$R^N(a^*,b_{\gamma_1}) \leftrightarrow
 R^N(a^*,b_{\gamma_2})$.

% If not, suppose $R^N(a,b_{\gamma_1})$ but not
%$R^N(a,b_{\gamma_2})$. Then $b_{\gamma_1} - b_{\gamma_2} \in A$ but
%$R(N,b_{\gamma_1} - b_{\gamma_2}) \cap P^{M_{\alpha^*}}_0 =
%\emptyset$, contradicting the choice of $\alpha^*$.

Continuing the proof of case ii) we define an ultrafilter  $\Dscr$ on
$\Pscr(P^{M_{\alpha^*}}_0)$ by $v \in \Dscr$ if for some (and hence
any) $b_\gamma$  {\em from our chosen subsequence} with $R(M,b_\gamma) \cap P^{M_{\alpha^*}}_0 = v$,
$R^N(a^*,b_\gamma)$. (This is an ultrafilter as each $u \subset
P^{M_{\alpha^*}}_0$ is $R(M,b_\gamma) \cap P^{M_{\alpha^*}}_0$  for some $\gamma$ by
requirement 4  of Construction~\ref{goalzfc+}.)

Now we show the coding of the elements of $\Dscr$ extends to the
entire original sequence.

 \begin{claim}\label{easy} For any $b \in P^{M}_1$, which is one of the original sequence of independent $b_\alpha$, if
  $v =R(M,b) \cap
P^{M_{\alpha^*}}_0$ and $v \in \Dscr$ then $N\models R(a^*,b)$.
\end{claim}

Proof. We can choose  $\beta, \beta_1$   so that $\alpha^*< \beta
< \lambda$, $b \in P^{M_{\beta}}_1$ and  $\beta_1 > \beta$
such that
$v_{\beta_1} = v$. % , v_{\beta_2}  = v?$.
Now $\check b = b \triangle b_{\beta_1} \in P_1^M$ and $R(M,\check b)
\cap P^{M_{\alpha^*}}_0 = \emptyset$.
Note that since $\langle b_\beta\colon \beta < \lambda\rangle$ is
 independent from  $P^{M_{\alpha^*}}_1$ over $P^M_4$ in $P^M_1$, in
 particular $b$ and $b_{\beta_1}$ are independent  so the singleton
  $b \triangle b_{\beta_1}$ is independent from  $P^{M_{\alpha^*}}_1$ over $P^M_4$ in $P^M_1$. So by the choice of $\alpha_*$,
$N \models \neg R(a^*,\check b)$.  So, $N \models \neg R(a^*,b)$ if
and only $N \models \neg R(a_*,b_{\beta_1})$.  But, we have $v \in
\Dscr$ and $R(M,b_{\beta_1}) \cap P_0^{M_{\alpha_*}} = v$, so $N
\models R(a^*,b_{\beta_1})$ and thus $N \models R(a_*,b)$ as
required.  $\qed_{\ref{easy}}$

\medskip
There is no
 $\aleph_1$-complete ultrafilter on $\Pscr(P^{M_{\alpha^*}}_0)$ since $|P^{M_{\alpha^*}}_0| < \lambda$
is not measurable. So there are $\langle w_n \subseteq
P^{M_{\alpha^*}}_0\colon n< \omega\rangle$, each in $\Dscr$, that are
decreasing and intersect in $\emptyset$.
Now we can find $\delta^*> \alpha^*$ such that $\delta^* \in S_A \cap
C$, the universe of $M_{\delta^*}$ is $\delta^*$,
$A^{}_{\delta^*} \cap \delta^* = A \cap \delta^*$, %= A_{\delta_v}$.
and there is an increasing sequence $\langle\gamma^{\delta^*}_n: n< \omega\rangle$
with limit at most $\delta^*$ and each $\gamma^{\delta^*}_n\not \in S$.
Further, by requirement 4 on the construction,  we can choose ${\gamma^{\delta^*}_n}$ so that
 $b_{{\gamma^{\delta^*}_n}}$ (another subsequence of the orginal sequence) satisfies
 $a \in R(M_{\gamma_n},b_{\gamma_n})$ if and only if $a \in w_n$,
%$R(M,b_{\gamma^{\delta^*}_n})\cap M_{\alpha^*} = w_n$,
$b_{\gamma^{\delta^*}_n} \in M_{\gamma^{\delta^*}_{n+1}}$, and the
sequence $\{b_{\gamma^{\delta^*}_{n}}\}$ is independent from
$P^{M_{\delta^*}}_1$ over $P^M_4$.
Since the $w_n$ are decreasing with empty intersection, no $a\in M_{\alpha^*}$ is in more than finitely many of the $w_n$.
%Further,
%for any $a \in M$ and any $n$,
%%and $m>n$,  $a \in R(M,b_{\gamma_n})$ if and only if
%$a \in R(M_{\gamma_n},b_{\gamma_n})$ if and only if $a \in w_n$.
Thus, Definition \ref{gooddef}~\ref{empint} is satisfied.

\medskip
So by clause 5) of the Requirements~\ref{goalzfc+}, there is a {\em labeled} $A$-good sequence $ \hat
     b_{\delta^*,n}$ for
  $M_{\delta^*+1}$, $\hat c_\delta^* \in
     P^{M_{\delta^*+1}}_2$ such that for each $n$,
     $F^{M_{\delta^*+1}}_n(\hat c_\delta^* ) = \hat
     b_{\delta^*,n}$. And by clause 8 of Definition~\ref{defK1}, this contradicts
     Claim~\ref{easy}; the intersection of $R(N,F^N_n(c))$ for $n< \omega$ must be empty but it contains $a^*$.
 So we
finish case ii) and thus Lemma~\ref{zfc+thm}. $\qed_{\ref{zfc+thm}}$
}
\end{construction} % $\qed_{Details \ref{verzfc+}}$
\medskip

%\sidebar{Aug 26: Problem  A2-- in previous paragraph and in March 6, 2017 version page 14 $\bigoplus_6$.  ($w_n$ here is $v_n$ in SS)
%
%$b_{\gamma_n}$ means $R(M_{\gamma_n}, b_{\gamma_n})$ in SS
%
%i) $b_{\gamma_n} \cap \gamma_n = w_n$
%
%and ii) the $b_{\gamma_n}$ are independent.
%
%
%Since the $w_n$ are descending this is contradictory????
%
%When $b_{\gamma_{n+1}}$ is chosen $w_n \subset R(M_{\gamma_{n+1}}, b_{\gamma_n})$ and
%$R(M_{\gamma_{n+1}},b_{\gamma_{n+1} }$ is set equal to $w_{n+1} \subset w_n$.
%But $R(M_{\gamma_{n+1},b_{\gamma_n}})$ and
%supposed to be independent since $b_{\gamma_{n+1}}$ and $b_{\gamma_{n}}$ are.$R(M_{\gamma_{n+1}},b_{\gamma_{n+1} })$ are
%}

%But the existence of such a sequence violates Claim~\ref{cl3} s

%\sidebar{Here is the issue. Since $\widehat \bK$ is claimed to have amalgamation with respect to $\subseteq_eq$ in all powers, we cannot show $\bK$ has
%no extension in $\widehat \bK$.  So where do we use that the extension is in $\Krich$?
%
%In fact claim 1.19 of oct 16 says it is maximal for $\K_0$ which I take as a typo for $\bK_1$ which is the earlier notation for $\Krich$ .
%}

 %Clearly,
%$\neg P_4(b_0)$. Since by Lemma~\ref{K2elem}.2,
%$M{\prec_{\infty,\omega}} N$, $b_0$ also generates a principal ideal
% in $N$. Since  $R^N(a^*,b_0)$, in $N$, $b_0 \leq G_1(a_*)$ so $b_0 \in P^N_4$. This is a contradiction.

\begin{remark}\label{context}  %The argument for the main claim does not explicitly rely on the
%claim that $M \in \bK_2$ rather than merely $\bK_1$.
 %But unless we
%knew $M$ were free, we would not be able to apply
%Corollary~\ref{getfreeext} to obtain an extension of $M$ in $\bK_2$
%and thus the actual model.
{\rm
In the construction we showed for limit $\delta$ that $M_\delta$ is $\bK_1$-free using $S$ does not reflect if $\delta \not \in S$ and that
$\cf(\delta) =\omega$ for $\delta\in S$.
We have no such    tools to  show the $P_0$-maximal model, $M =M_\lambda$
built in Theorem~\ref{mainclaim}    is $\bK_1$-free. In fact, by the
contrapositive of Corollary~\ref{getfreeext} the final  $P_0$-maximal model,
which might be $M$, is {\em not} $\bK_1$-free.

Note that every subset of $M$ with cardinality $< \lambda$ is contained in a $\bK_1$-free substructure; this fails in the ZFC proof \cite{BaldwinShelahhanfmaxzfc} of maximal models of $\bK_2$ cofinal in a measurable.
}
\end{remark}

\medskip

Recall that a $P_0$-maximal model in a class $\bK$ is one that cannot be extended in $\bK$ without
extending $P_0$. While a maximal model has no extension $\bK$. We have constructed a $P_0$-maximal model in $\bK_2$; we show that it has a
$\bK_2$-maximal extension that is only slightly larger.

\begin{corollary}\label{mainclaim} Under the hypotheses of Theorem~\ref{zfc+thm},
there is a maximal model of $\bK_2$ of cardinality at most
$2^\lambda$.
\end{corollary}

Proof.  Fix a $P_0$-maximal model $N_0$ of cardinality $\lambda$ from
Theorem~\ref{zfc+thm}. Build for as long as possible a continuous
$\subseteq$-increasing chain of $N_\alpha \in \bK_2$ such that each
$P_1^{N_\alpha} \subsetneq P_1^{N_{\alpha +1}}$. But, necessarily, $P_0^{N_\alpha} = P_0^{N_{\alpha +1}}$. Recall that by Lemma~\ref{nonoise}.1
the relation $R$ is injective. So, each $|P_1^{N_\alpha}| \leq
2^{|P_0^{N_0}|}=2^\lambda$. So this construction must stop and the
final, maximal in $\bK_2$, model has cardinality at most
$2^\lambda$. $\qed_{\ref{mainclaim}}$

\section{Hanf Number for Existence}\label{bkl}

As mentioned in the introduction, we improved in \cite{BKL} Hjorth's
result \cite{Hjorthchar} by exhibiting for each $n< \omega$ a
complete sentence $\psi_n$ such that $\psi_n$ characterizes
$\aleph_n$.  This improvement is achieved by combining the combinatorial
idea of Laskowski-Shelah in \cite{LaskowskiShelahatom} with a new
notion of $n$-dimensional amalgamation. We explain the main
definition and theorem here (as in the Tehran lectures) and refer to
\cite{BKL} for the proofs.
The combinatorial fact is:

\begin{fact}\cite{LaskowskiShelahatom}\label{comb1}  For every $k\in\omega$, if $\cl$ is a locally finite
closure relation on a set $X$ of size $\aleph_k$, then there is an
independent subset of size $k+1$.
\end{fact}

 Fix a vocabulary $\tau_r$ with infinitely many $r$-ary relations $R_n$ and
infinitely many $r+1$-ary functions $f_n$.
We consider the class $\bK^r_0$ of finite $\tau_r$-structures
(including the empty structure) that satisfy the following three
conditions; closure just means subalgebra closure with respect to the
functions.

\begin{itemize}
\item The relations $\{R_n\colon n\in\omega\}$ partition the
    $(r+1)$-tuples;
\item  For every $(r+1)$-tuple $\abar=(a_0,\dots,a_r)$, if
    $R_n(\abar)$ holds, then $f_m(\abar)=a_0$ for every $m\ge n$;
\item  There is no independent subset of size $r+2$.
\end{itemize}

It is easy to see from Fact~\ref{comb1} that every model in
$\aleph_r$ is maximal. The main effort is to show there is a complete
sentence $\phi_r$ satisfying those conditions which has model in
$\aleph_{r}$.  For this we introduce a notion patterned on
excellence\footnote{Shelah's theory of excellence concerns unique
free disjoint amalgamations of infinite structures in $\omega$-stable
classes of models of complete sentences in $L_{\omega_1,\omega}$.}
but weaker.
% This extension was necessary as the proof in the
%excellence context that $\aleph_{\alpha,n}$ amalgamation implies
%$\aleph_{\alpha+1,n-1}$ amalgamation {\em depends essentially on
%Fodor's lemma;} a new framework was needed to pass from finite to
%countable.
We pass from a class $\bK_0^r$ of, now, locally finite
structures to the associated class $\widehat \bK$ as in
Definition~\ref{defclass}.

\begin{definition}  For $k\ge 1$, a {\em $k$-configuration} is a sequence $\Mbar=\langle M_i:i<k\rangle$ of models (not isomorphism types)
 from $\bK$.
We say {\em $\Mbar$ has power $\lambda$} if $\|\bigcup_{i<k}
M_i\|=\lambda$. An {\em extension} of $\Mbar$ is any $N\in\bK$ such
that every $M_i$ is a substructure of $N$.
\end{definition}

Informally, $(\lambda,k)$-disjoint amalgamation holds when for any
sequence of $k$ models, at least one with $\lambda$ elements, there is
common extension, which properly extends each model in the sequence.  Crucially, there is no prior
assumption of a universal model.  Here is the precise formulation.

\begin{definition}
Fix a cardinal $\lambda=\aleph_{\alpha}$ for $\alpha\ge -1$.  We
define the notion of a class $(\bK,\le)$ having {\em
$(\lambda,k)$-disjoint amalgamation} in two steps:
\begin{enumerate}
\item  $(\bK,\le)$ has $(\lambda,0)$-disjoint amalgamation if
    there is $N\in\bK$ of power $\lambda$;
\item  For $k\ge 1$, $(\bK,\le)$ has $(\leq \lambda,k)$-disjoint
    amalgamation if it has $(\lambda,0)$-disjoint amalgamation
    and every $k$-configuration $\Mbar$ of cardinality
    $\leq\lambda$ has an extension $N\in\bK$ such that every
    $M_i$ is a proper substructure of $N$.
\end{enumerate}
For $\lambda\ge\aleph_0$, we define $(<\lambda,k)$-disjoint
amalgamation by:  has $(\leq \mu,k)$-disjoint amalgamation for each
$\mu<\lambda$.
\end{definition}

Whether or not a given $k$-configuration $\Mbar$ has an
    extension depends on more than the sequence of isomorphism
    types of the constituent $M_i$'s, as the pattern of
    intersections is relevant as well. For example, when (as here) strong substructure
    is just substructure), a
    2-configuration $\<M_0,M_1\>$ with neither contained in the
    other has an extension   if and only if the triple of
    structures $\<M_0\cap M_1,M_0,M_1\>$ has an extension
    amalgamating them disjointly. Thus we abuse notation a bit and write $(<\lambda,2)$ amalgamation
    for both the notion defined here and the one in Definition~\ref{dap2}. But there is no existing analog
    of our disjoint $(<\lambda, k)$-amalgamation for $k>2$.

Now we modify a theme familiar from the theory of excellence. If the
cardinality increases by one,   the number of models that can be
amalgamated drops by one. In  Shelah's context \cite{Shaecbook}
(chapter 21 of \cite{Baldwincatmon}) there is a reliance  on Fodor's
lemma to obtain compatible filtrations of the models in $\kappa^+$ to
prove the version of Proposition~\ref{trans}.  A very different
approach was needed to go from the finite to the countable.  Instead
of the $k$th level concerning finding an embedding into an upper
corner for a given $2^{k-1}$ vertices of a $k$-cube, we consider
actual containment for $k$-models and do not worry about their
intersections.

\begin{lemma} [Proposition 2.20 of \cite{BKL}] \label{trans}
 Fix a locally finite $(\bK,\le)$  with JEP.  For all
cardinals $\lambda\ge\aleph_0$ and for all $k\in\omega$, if $\bK$ has
$(<\lambda,k+1)$-disjoint amalgamation, then it also has
$(\le\lambda,k)$-disjoint amalgamation.
\end{lemma}

Together, these propositions yield 1)-3)  of the next result.  Recall from Definition~\ref{dap2}, that by $2$-amalgamation, we mean the usual notion that allows identifications. We say $2$-amalgamation is {\em trivially} true in a cardinal $\kappa$ if all models in $\kappa$ are maximal.

\begin{theorem} [Theorem 3.2.4 of \cite{BKL}] \label{sumup}
For every $r \ge 1$, the class ${\bf At^r}$ satisfies:
\begin{enumerate}

\item there is a  model of size $\aleph_r$, but no larger models;

\item every model of size $\aleph_r$ is maximal, and so
     $2$-amalgamation is trivially true in $\aleph_r$;

\item disjoint $2$-amalgamation holds up to
    $\aleph_{r-2}$;
%\item The situation for atomic models (or $\hat \bK^{r}$)
%    $\aleph_{r-1}$ is unclear
%    \begin{enumerate}
\item $2$-amalgamation fails in $\aleph_{r-1}$.

\item Each of the classes $\hat {\bK^r}$ and ${\bf At}^r$ have
    $2^{\aleph_s}$ models in $\aleph_s$ for $1 \leq s \leq r$. In
    addition,  $\hat {\bK^r}$ has $2^{\aleph_0}$ models in
    $\aleph_0$.
\end{enumerate}

\end{theorem}

Parts 4) and 5) require a further refinement of the notion of disjoint amalgamation.

\begin{definition}  \label{frugal}
Given a cardinal $\lambda$ and $k\in\omega$, we say that
    $\bK$ has {\em frugal $(\leq\lambda,k)$-disjoint
    amalgamation} if it has $(\leq\lambda,k)$-disjoint
    amalgamation and, when $k\ge 2$, every $k$-configuration
    $\<M_i:i<k\>$ of cardinality $\leq \lambda$ has an extension $N\in\bK$ with universe
    $\bigcup_{i<k} M_i$.

\end{definition}

Thus the domain of a  frugal amalgamation is just the union of the models amalgamated. It is easy to see that
this property holds for the example in \cite{BKL}. It is essential for the intricate constructions to verify the last two parts of Theorem~\ref{sumup}
and for the work in \cite{BKSoul, BSoul}.

The {\em finite amalgamation spectrum} of an abstract elementary
class $\bK$  with $LS(\bK)= \aleph_0$ is the set $X_{\bK}$ of $n<\omega$
such that    $\bK$ satisfies
amalgamation\footnote{We say amalgamation holds in $\kappa$ in the
trivial special case when all models in $\kappa$ are maximal.
 %to distinguish that situation
%from there simply being no models in $\kappa^+$.
We say amalgamation fails in $\kappa$ if there are no models to
amalgamate.} in $\aleph_n$.  There are many examples\footnote{Kueker,
as reported in \cite{Malitz}, gave the first example of a complete
sentence failing amalgamation in $\aleph_0$.} where the finite
amalgamation spectrum of a complete sentence  of
$L_{\omega_1,\omega}$ is either $\emptyset$ or $\omega$.

Theorem~\ref{sumup}  gave the
first example of such a sentence with a non-trivial spectrum: for each $1\le r < \omega$
amalagmation holds up to $\aleph_{r-2}$, but fails in $\aleph_{r-1}$.  It holds (trivially) in
$\aleph_{r}$ (since all models are maximal); there is no model in
$\aleph_{r+1}$.

This result leaves open whether the property, AP in $\lambda$, can be
true or false in various patterns as $\lambda$ increases? Is there
even an AEC (and more interestingly a complete sentence of
$L_{\omega_1,\omega}$) and
 cardinals $\kappa<\lambda$ such that amalgamation holds non-trivially in both $\kappa$ and $\lambda$
 but fails at some cardinal between them?

 Relying on the construction in \cite{BKL}, Baldwin and
Souldatos \cite{BSoul} show there exist {\em complete} sentences of
$L_{\omega_1,\omega}$ that variously have  maximal models a) in two
successive cardinals, b)  in $\kappa$ and $\kappa^{\omega}$ and c) in
countably many cardinals.  In each case all maximal models of the
sentence have cardinality less than $\aleph_{\omega_1}$.  That proof
includes an intricate construction of a complete sentence that has a
model in each successor cardinal $\kappa^+$ with a definable subset
of power $\kappa$.  The \cite{BSoul} result is distinguished from the
one here in several ways. It constructs maximal models in designated
cardinals rather than an initial segment.  The crucial amalgamation
properties are quite different. The example in \cite{BKL} satisfies
$(<\lambda,2)$ amalgamation in all cardinals.

%\sidebar{FIX} In \cite{BSoul} The most delicate argument in \cite{BS}
%shows that one can amalgamate  a model of $\widehat \bK$ of any
%cardinality with an arbitrary finite model and thus achieve richness.
%\bibliography{C:/Users/jbald/texmf/bibtex/bib/local/ssgroups}
%\bibliographystyle{alpha}
%
%\end{document}
%\bibliography{C:/Users/jbald/texmf/bibtex/bib/local/ssgroups}
%\bibliographystyle{alpha}
%
%\end{document}

\end{document}